\documentclass[leqno,12pt]{article}
\usepackage{}

\usepackage{amssymb}
\usepackage{euscript}
\usepackage[dvips]{graphicx}
\usepackage{flafter}
\usepackage{pstricks}
\usepackage[latin1]{inputenc}
\usepackage{amsmath}
\usepackage{amsfonts}

\makeatletter

\newcommand{\Rmnum}[1]{\expandafter\@slowromancap\romannumeral #1@}
\makeatother

\usepackage{verbatim}

\usepackage{mathrsfs} 

\evensidemargin\oddsidemargin

\begin{document}

\baselineskip=18pt
\setcounter{page}{1}

\renewcommand{\theequation}{\thesection.\arabic{equation}}
\newtheorem{theorem}{Theorem}[section]
\newtheorem{lemma}[theorem]{Lemma}
\newtheorem{proposition}[theorem]{Proposition}
\newtheorem{corollary}[theorem]{Corollary}
\newtheorem{remark}[theorem]{Remark}
\newtheorem{fact}[theorem]{Fact}
\newtheorem{problem}[theorem]{Problem}
\newtheorem{example}[theorem]{Example}
\newtheorem{question}[theorem]{Question}
\newtheorem{conjecture}[theorem]{Conjecture}

\newcommand{\eqnsection}{
\renewcommand{\theequation}{\thesection.\arabic{equation}}
    \makeatletter
    \csname  @addtoreset\endcsname{equation}{section}
    \makeatother}
\eqnsection

\def\r{{\mathbb R}}
\def\e{{\mathbb E}}
\def\p{{\mathbb P}}
\def\P{{\bf P}}
\def\E{{\bf E}}
\def\Q{{\bf Q}}
\def\z{{\mathbb Z}}
\def\N{{\mathbb N}}
\def\T{{\mathbb T}}
\def\G{{\mathbb G}}
\def\L{{\mathbb L}}

\def\deg{\chi}

\def\ee{\mathrm{e}}
\def\d{\, \mathrm{d}}
\def\S{\mathscr{S}}



\vglue50pt

\centerline{\large\bf Waiting times for particles in a branching Brownian motion }
\centerline{\large\bf to reach the rightmost position }

\bigskip
\bigskip

\centerline{by}

\medskip

\centerline{Xinxin Chen}

\medskip

\centerline{\it Universit\'e Paris VI}

\bigskip
\bigskip
\bigskip

{\leftskip=2truecm \rightskip=2truecm \baselineskip=15pt \small

\noindent{\slshape\bfseries Summary.} It has been proved by Lalley and Sellke \cite{Lalley1987} that every particle born in a branching Brownian motion has a descendant reaching the rightmost position at some future time. The main goal of the present paper is to estimate asymptotically as $s$ goes to infinity, the first time that every particle alive at the time $s$ has a descendant reaching the rightmost position.

\bigskip

\bigskip


} 

\bigskip
\bigskip

\section{Introduction}
   \label{s:intro}
\subsection{The model}
We consider a branching Brownian motion (BBM) on the real line $\r$, which evolves as follows. Starting at time $t=0$, one particle located at 0, called the root, moves like a standard Brownian motion until an independent exponentially distributed time with parameter 1. At this time it splits into two particles, which, relative to their birth time and position, behave like independent copies of their parent, thus moving like Brownian motions and branching at rate 1 into two copies of themselves. Let $\mathcal{N}(t)$ denote the set of all particles alive at time $t$ and let $N(t):=\#\mathcal{N}(t)$. For any $v\in\mathcal{N}(t)$  let $X_v(t)$ be the position of $v$ at time $t$; and for any $s<t$, let $X_v(s)$ be the position of the unique ancestor of $v$ that was alive at time $s$. We define
$$R(t):=\max_{u\in\mathcal{N}(t)}X_u(t)\text{ and }L(t):=\min_{u\in\mathcal{N}(t)}X_u(t),$$
\noindent which stand for the rightmost and leftmost positions, respectively.

The positions of the extremal particles of a BBM, $R(t)$, have been much studied both analytically and probabilistically. Kolmogorov et al. \cite{Kolmogorov1937} proved that $R(t)/t$ converges almost surely to $\sqrt{2}$. Bramson \cite{Bramson1978} \cite{Bramson1983} showed that $R(t)-\sqrt{2}t+(3/2\sqrt{2})\log t$ converges in law. These results hold as well for a wide class of branching random walks under mild conditions: see for example Biggins \cite{Biggins1976}, Addario-Berry and Reed \cite{Reed2009}, Hu and Shi \cite{Hu2009}, A\"{\i}d\'ekon \cite{Aidekon2011}. In particular, we state the following fact, which is first given by Hu and Shi \cite{Hu2009} for branching random walks, and is recently proved by Roberts \cite{Roberts2011}:
\begin{eqnarray}
\liminf_{t\rightarrow\infty}\frac{R(t)-\sqrt{2}t}{\log t} &=& -\frac{3}{2\sqrt{2}}\quad\text{  almost surely};\label{liminf}\\
\limsup_{t\rightarrow\infty}\frac{R(t)-\sqrt{2}t}{\log t} &=& -\frac{1}{2\sqrt{2}}\quad\text{  almost surely}.\label{limsup}
\end{eqnarray}

In \cite{Lalley1987}, Lalley and Sellke showed the following interesting property: every particle born in a BBM has a descendant reaching the rightmost position at some future time.  Such a particle was thought of having a prominent descendant ``in the lead" at this time. This property is in agreement with the branching-selection particle systems investigated in \cite{Brunet1997}, \cite{Brunet1999} and \cite{Berard2008}.  These papers bring out the fact that the extremal positions of a branching system on the line cannot always be occupied by the descendants of some ``elite" particles.

In the present work, we give some quantitative understanding of this behavior, and precisely speaking, about how long we have to wait so that every particle alive at time $s$ has a descendent that has occupied the rightmost position.
\subsection{The main problem}
Let us make an analytic presentation for our problem. For any $s>0$ and each particle $u\in\mathcal{N}(s)$, the shifted subtree generated by $u$ is
\begin{equation}\label{shiftedsubtree}
\mathcal{N}^u(t):=\Big\{v\in\mathcal{N}(t+s), u\leq v\Big\},\quad \forall t\geq 0,
\end{equation}
\noindent where $u\leq v$ indicates that $v$ is a descendant of $u$ or is $u$ itself.
Further, for any $v\in\mathcal{N}^u(t)$, let
\begin{equation}\label{shiftedpositions}
X^u_v(t):=X_v(t+s)-X_u(s),
\end{equation}
\noindent be its shifted position.  We set $R^u(t):=\max_{v\in\mathcal{N}^u(t)}X_v^u(t)$ and $L^u(t):=\min_{v\in\mathcal{N}^u(t)}X^u_v(t)$. Moreover, Let $\Big\{\mathcal{F}_t; t\geq 0\Big\}$ be the natural filtration of the branching Brownian motion. The branching property implies that, given $\mathcal{F}_s$, $\{R^u(\cdot); u\in\mathcal{N}(s)\}$ are independent copies of $R(\cdot)$. Moreover, we denote by $\mathcal{F}_\infty^u$ the sigma-field generated by the shifted subtree started from the time $s$ rooted at $u$.

For every $u\in\mathcal{N}(s)$, let
\begin{equation}\label{defofone}
\tau_u:=\inf\{t>0: R(t+s)=X_u(s)+R^u(t)\}.
\end{equation}
 The random variable $\tau_u$ stands for the first time that started from time $s$, the particle $u$ has a descendant reaching the rightmost position in the system. It is the object in which we are interested. We define
\begin{equation}
\Theta_s:=\max_{u\in\mathcal{N}(s)}\tau_u,
\end{equation}
which represents the first time when every particle in $\mathcal{N}(s)$ has had a descendant occupying the rightmost position.

According to Lalley and Sellke \cite{Lalley1987}, for any $s>0$, $\mathbb{P}[\Theta_s<\infty]=1.$
\noindent Since $\Theta_s\rightarrow\infty$ almost surely as $s\rightarrow\infty,$ we intend to determine the rate at which $\Theta_s$ increases to infinity.
\subsection{The main results}
To estimate $\Theta_s=\max_{u\in\mathcal{N}(s)}\tau_u$, an intuitive idea consists in saying that, the further a particle is away from the rightmost one, the longer it has to wait for a descendant to be located on the rightmost position. We thus first focus on the leftmost particle.  Let $\ell(s)$ be the leftmost particle alive at time $s$. By (\ref{defofone}), $\tau_{\ell(s)}$ is defined as the shortest time needed for $\ell(s)$ to wait to have a descendant occupying the rightmost position.
\begin{theorem}\label{leftmost particle}
The following convergence holds almost surely
\begin{equation}\label{leftmost particleeq}
\lim_{s\rightarrow\infty}\frac{\log\tau_{\ell(s)}}{s}=4.
\end{equation}
\end{theorem}

However, the leftmost particle is not the one who ``drags the feet" of the whole population $\mathcal{N}(s)$. By considering the positions of all particles alive at time $s$, as well as their evolutions, we obtain our main result as follows.
\begin{theorem}\label{mainconclusion}
The following convergence holds almost surely
\begin{equation}\label{mainconclusioneq}
\lim_{s\rightarrow\infty}\frac{\log\Theta_s}{s}=2+2\sqrt{2}>4.
\end{equation}
\end{theorem}

\begin{remark}
The proof of the theorems will reveal that the largest $\tau_u$ for $u\in\mathcal{N}(s)$ is achieved by some particle located at a position around $-(2-\sqrt{2})s$ which does not split until time $s+\frac{1}{\sqrt{2}}s$ and moves towards to the left as far as possible.
\end{remark}

The rest of this paper is organized as follows. Section 2 is devoted to discussing the behaviors of the extremal position $R(\cdot)$, which leads to two propositions. In Section 3, we consider the case of two independent branching Brownian motions and state another proposition. We prove Theorem \ref{leftmost particle} in Section 4 by means of these propositions. Finally, in Section 5, we prove Theorem \ref{mainconclusion}.

\section{The behavior of the rightmost position}
   \label{s:rightmost}
$\phantom{aob}$We recall Proposition 3 in Bramson's work \cite{Bramson1978}. It is shown that for all $0\leq y\leq t^{1/2}$ and $t\geq 2$, there exists a positive constant $c$ which is independent of $t$ and $y$, such that
\begin{equation}
\mathbb{P}\bigg[R(t)>m(t)+y\bigg]\leq c(1+y)^2\exp(-\sqrt{2}y),
\end{equation}
where
$$m(t):=\sqrt{2}t-\frac{3}{2\sqrt{2}}\log t.$$
Therefore, with $c_1:=c+1$, we get the following inequality, which will be applied several times in our arguments.

\begin{fact}[Bramson \cite{Bramson1978}]  For any $t\geq 2$ and $y\leq \sqrt{t}$,
\begin{equation}\label{ineq}
\mathbb{P}\bigg[R(t)>m(t)+y\bigg]\leq c_1\Big(1+y_+\Big)^2e^{-\sqrt{2}y},
\end{equation}
with $y_+:=\max\{y, 0\}$.
\end{fact}

Let $(B_s; s\geq0)$ be a standard Brownian motion on $\mathbb{R}$. We state the following lemma, which can be found in several papers (e.g. \cite{Lyons1995} \cite{Harris2011}). It is also of frequent use.
\begin{lemma}[many-to-one]
For any measurable function $F$ and each $t>0$,
\begin{equation}
\mathbb{E}\bigg[\sum_{u\in\mathcal{N}(t)}F(X_{u}(s), s\in[0,t])\bigg]=e^t\mathbb{E}\bigg[F(B_s, s\in[0,t])\bigg],
\end{equation}
where, for each $u\in\mathcal{N}(t)$ and $s\in[0,t]$, $X_{u}(s)$ denotes the position, at time $s$, of the
ancestor of $u$.
\end{lemma}

Let us present the following inequality as well, which is Equation (57) in Bramson \cite{Bramson1978}.
\begin{fact}For any $s\geq 1$ and any $a>0$,
\begin{equation*}
\Big(1-\frac{s}{a^2}\Big)\sqrt{\frac{s}{2\pi}}a^{-1}\exp\Big(-\frac{a^2}{2s}\Big)\leq\mathbb{P}[B_s\geq a ]\leq \sqrt{\frac{s}{2\pi}}a^{-1}\exp\Big(-\frac{a^2}{2s}\Big).
\end{equation*}
\end{fact}
It immediately follows that
\begin{equation}\label{evaluationofBM}
\mathbb{P}[B_s\leq -a]=\mathbb{P}[B_s\geq a ]\leq\frac{\sqrt{s}}{a}\exp\Big(-\frac{a^2}{2s}\Big).
\end{equation}
Moreover, if $a=\alpha s$ with some constant $\alpha>0$, we have
\begin{equation}\label{convergenceofBM}
\mathbb{P}[B_s\leq -\alpha s]=\mathbb{P}[B_s\geq \alpha s ]=\exp\Big\{-\Big(\frac{\alpha^2}{2}+o_s(1)\Big)s\Big\},
\end{equation}
where $o_s(1)\rightarrow0$ as $s$ goes to infinity.

We define, for any $y>0$,
$$T(y):=\inf\Big\{t\geq 1; R(t)-m(t)>y\Big\}.$$
Because of (\ref{limsup}), one immediately sees that $\mathbb{P}\Big[T(y)<\infty\Big]=1$ for any $y>0$. Moreover, $T(y)\uparrow\infty$ almost surely as $y\uparrow\infty$.

\begin{proposition}\label{stoppingtime}
The following convergence holds almost surely
\begin{equation}\label{stoppingtime}
\lim_{y\rightarrow\infty}\frac{\log T(y)}{y}=\sqrt{2}.
\end{equation}
\end{proposition}

\textbf{Proof:} First, we prove the lower bound.

Let $2\leq y\leq \sqrt{t}$, and set
\begin{equation*}
\Lambda:=\mathbb{E}\Bigg[\int_1^{t+1}\mathbf{1}_{(R(s)>m(s)+y-1)}ds\Bigg].
\end{equation*}
\noindent Clearly,  $\Lambda=\int_{1}^{t+1} \mathbb{P}\Big[R(s)>m(s)+y-1\Big] ds$. Hence,
\begin{eqnarray*}
\Lambda &= & \int_{1}^{y^2} \mathbb{P}\Big[R(s)>m(s)+y-1\Big] ds+\int_{y^2}^{t+1} \mathbb{P}\Big[R(s)>m(s)+y-1\Big] ds\\
                    & \leq & \int_{1}^{y^2} \mathbb{E}\bigg[\sum_{u\in\mathcal{N}(s)}1_{(X_u(s)>m(s)+y-1)}\bigg] ds+\int_{y^2}^{t+1} \mathbb{P}\Big[R(s)>m(s)+y-1\Big] ds.
\end{eqnarray*}
By the many-to-one lemma and by (\ref{evaluationofBM}),
\begin{eqnarray}\label{smalltimepre}
\int_{1}^{y^2} \mathbb{E}\bigg[\sum_{u\in\mathcal{N}(s)}1_{(X_u(s)>m(s)+y-1)}\bigg] ds &=& \int_{1}^{y^2} e^s\mathbb{P}\Big[B_s>m(s)+y-1\Big]ds\nonumber\\
                                                                                                                                            &\leq&\int_1^{y^2} \frac{ \sqrt{s}}{m(s)+y-1} \exp\Big\{\frac{-(m(s)+y-1)^2}{ 2s}+s\Big\}ds.
                                                                                                                                            \end{eqnarray}
Note that for $m(s)=\sqrt{2}s-\frac{3}{2\sqrt{2}}\log s$ with $s\in[1,y^2]$, the inequalities
\begin{equation}
m(s)+y-1\geq \sqrt{2}s\text{ and } \exp(\frac{-(m(s)+y-1)^2}{ 2s}+s)\leq s^{3/2}e^{-\sqrt{2}(y-1)}
\end{equation}
hold. Plugging them into the integration of (\ref{smalltimepre}) yields that
\begin{equation}\label{smalltime}
\int_{1}^{y^2} \mathbb{E}\bigg[\sum_{u\in\mathcal{N}(s)}1_{(X_u(s)>m(s)+y-1)}\bigg] ds\leq \int_1^{y^2} \frac{s^2}{\sqrt{2}s}e^{-\sqrt{2}(y-1)}ds\leq c_2 y^4e^{-\sqrt{2}y},
\end{equation}
which is then bounded by $c_2t y^2 e^{-\sqrt{2}y}$ as $y\leq\sqrt{t}$.
Meanwhile, by the inequality (\ref{ineq}),
\begin{equation}\label{bigtime}
\int_{y^2}^{t+1} \mathbb{P}\Big[R(s)>m(s)+y-1\Big] ds\leq \int_{y^2}^{t+1} c_1 y^2e^{-\sqrt{2}y+\sqrt{2}} ds\leq c_1 t y^2 e^{-\sqrt{2}y+\sqrt{2}}.
\end{equation}
Combining (\ref{smalltime}) with (\ref{bigtime}), we have
\begin{equation}\label{smallpart}
\Lambda\leq c_3 t y^2e^{-\sqrt{2}y}.
\end{equation}
with $c_3>0$ a constant independent of $(y,t)$.

On the other hand,
\begin{eqnarray*}
\Lambda&\geq&\mathbb{E}\Bigg[\int_1^{t+1}\mathbf{1}_{(R(s)>m(s)+y-1)}ds;\, T(y)\leq t\Bigg]\\
                  & = &\int_{1}^t \mathbb{P}\Big[T(y)\in dr\Big]\mathbb{E}\Bigg[\int_1^{t+1}\mathbf{1}_{(R(s)>m(s)+y-1)}ds\Big\vert T(y)=r\Bigg].
\end{eqnarray*}
Conditionally on the event $\{T(y)=r\leq t\}$, the rightmost particle in $\mathcal{N}(r)$, denoted by $\omega$, is located at $m(r)+y$. Started from the time $r$, $\omega$ moves according to a Brownian motion and splits into two after an exponential time.
By ignoring its branches, we observe that $\Big\{R(s+r)>[m(r)+y]+[\sqrt{2}s-1]\geq m(s+r)+y-1\Big\}$ is satisfied as long as the Brownian motion realized by $\omega$ keeps lying above $\sqrt{2}s-1$. Hence, given $\{T(y)=r\leq t\}$,
\begin{equation*}
\int_1^{t+1}\mathbf{1}_{\big(R(s)>m(s)+y-1\big)}ds\geq_{st}\int_{0}^{t+1-r}\mathbf{1}_{\big(B_s>\sqrt{2}s-1\big)}ds\geq \min\{1, T^{(-\sqrt{2})}_{-1}\},
\end{equation*}
where $\geq_{st}$ denotes stochastic dominance and $T^{(-\sqrt{2})}_{-1}:=\inf\{t\geq 0; B_t< \sqrt{2}t-1\}.$

These arguments imply that
\begin{equation}
\begin{array}{cl}
\Lambda   & \geq  \int_{1}^t \mathbb{P}[T(y)\in dr] \mathbb{E}\Big[\min\{1,T^{(-\sqrt{2})}_{-1}\}\Big]\\
                  & =:  c_{13}\mathbb{P}[T(y)\leq t],
                  \end{array}
\end{equation}
\noindent where $c_{13}:=\mathbb{E}\Big[\min\{1, T^{(-\sqrt{2})}_{-1}\}\Big]\in(0,\infty)$. Compared with (\ref{smallpart}), this tells us that
\begin{equation}\label{lesspart}
\mathbb{P}\Big[T(y)\leq t\Big]\leq c_5  t y^2e^{-\sqrt{2}y}, \text{ for }2\leq y\leq \sqrt{t},
\end{equation}
where $c_5:=\frac{c_3}{c_{13}}\in(0,\infty)$.

Taking $t=e^{\sqrt{2}y(1-\delta)}$ with $\delta\in(0,1)$ yields that
$$\sum_{k=1}^\infty \mathbb{P}\Big[T(k)\leq e^{\sqrt{2}k(1-\delta)}\Big]<\infty.$$
According to the Borel-Cantelli lemma,
\begin{equation*}
\liminf_{y\rightarrow\infty}\frac{\log T(y)}{y}\geq \sqrt{2},\text{ almost surely, }
\end{equation*}
\noindent proving the lower bound in the proposition.

To prove the upper bound, we recall that
$$R^u(t)=\max\{X_v^u(t); v\in\mathcal{N}^u(t)\},\  u\in\mathcal{N}(s).$$
\noindent Obviously, $R^u(t);u\in\mathcal{N}(s)$ are i.i.d. given $\mathcal{F}_s$, and are distributed as $R(t)$.

We fix $a_y\in(0,y)$ and define the measurable events
\begin{eqnarray*}
\Sigma_1&:=&\bigg\{L(a_y)\geq -2 a_y; N(a_y)\geq \exp\Big(\frac{1}{2}a_y\Big)\bigg\},\\
\Sigma&:=&\Sigma_1\cap \Big\{T(y)>e^{\sqrt{2}y(1+\delta)}\Big\}.
\end{eqnarray*}
Then,
\begin{equation}\label{idea1}
\begin{array}{cl}
\mathbb{P}\bigg[T(y)>e^{\sqrt{2}y(1+\delta)}\bigg]&\leq \mathbb{P}\Big[\Sigma_1^c\Big]+\mathbb{P}\Big[\Sigma_1\cap \Big\{T(y)>e^{\sqrt{2}y(1+\delta)}\Big\}\Big]\\
&\leq  \mathbb{P}\bigg[N(a_y)\leq \exp\Big(\frac{1}{2}a_y\Big)\bigg]+\mathbb{P}\bigg[ L(a_y)\leq -2 a_y\bigg]+\mathbb{P}\Big[\Sigma\Big].
\end{array}
\end{equation}

We choose $a_y=\frac{\delta}{4+2\sqrt{2}}y=:\delta_1 y$ from now on to evaluate $\mathbb{P}[\Sigma]$. Since $\Sigma_1\in\mathcal{F}_{a_y}$, for $y$ large enough so that $2e^{\sqrt{2}y(1+\frac{1}{2}\delta)}\leq e^{\sqrt{2}y(1+\delta)}-a_y$, we have
\begin{eqnarray*}
\mathbb{P}\Big[\Sigma\big\vert \mathcal{F}_{a_y}\Big]&\leq & \mathbf{1}_{\Sigma_1}\prod_{u\in\mathcal{N}(a_y)} \mathbb{P}\Big[R^u(r)\leq m(a_y+r)+y-X_u(a_y), \forall r\leq e^{\sqrt{2}y(1+\delta)}-a_y\big\vert\mathcal{F}_{a_y}\Big]\\
 & \leq &\mathbb{P}\bigg[R(r)\leq m(r)+y+2a_y+\sqrt{2}a_y, \forall r\in \Big[ e^{\sqrt{2}y(1+\frac{1}{2}\delta)},2e^{\sqrt{2}y(1+\frac{1}{2}\delta)}\Big]\bigg]^{e^{a_y/2}}\\
 &\leq & \mathbb{P}\bigg[R(r)\leq m(r)+\frac{1}{\sqrt{2}}\log r, \forall r\in\Big[ e^{\sqrt{2}y(1+\frac{1}{2}\delta)},2e^{\sqrt{2}y(1+\frac{1}{2}\delta)}\Big]\bigg]^{e^{a_y/2}}.
\end{eqnarray*}

At this stage, it is convenient to recall the proof of Proposition 15 of Roberts \cite{Roberts2011}, saying that there exists a constant $c^\prime>0$ such that for $y$ large enough,
$$\mathbb{P}\bigg[ \exists r\in \Big[ e^{\sqrt{2}y(1+\frac{1}{2}\delta)},2e^{\sqrt{2}y(1+\frac{1}{2}\delta)}\Big]: R(r)\geq m(r)+\frac{1}{\sqrt{2}}\log r\bigg]>c^\prime>0.$$

Thus, $\mathbb{P}\Big[\Sigma\Big]\leq (1-c^\prime)^{e^{a_y/2}}\leq \exp(-c^\prime e^{\delta_1y/2})$.

It remains to estimate $\mathbb{P}\big[N(a_y)\leq \exp(\frac{1}{2}a_y)\big]$ and $\mathbb{P}\big[ L(a_y)\leq -2 a_y\big]$. On the one hand, the branching mechanism tells us that for any $s\geq 0$, $N(s)$ follows the geometric distribution with parameter $e^{-s}$ (for example, see Page 324 of \cite{Mckean1975}). It thus yields that $\mathbb{P}\big[N(a_y)\leq \exp(\frac{1}{2}a_y)\big]\leq e^{-\delta_1y/2}$.
On the other hand, as shown in Proposition 1 of Lalley and Sellke \cite{Lalley1989}, for any $\mu\geq \sqrt{2}$ and $s>0$,
\begin{equation}\label{feynmankac}
\mathbb{P}\Big[L(s)\leq -\mu s\Big]=\mathbb{P}\Big[R(s)\geq\mu s\Big]\leq \mu^{-1}(2\pi s)^{-1/2}\exp\bigg(-s\Big(\frac{\mu^2}{2}-1\Big)\bigg).
\end{equation}

Consequently, (\ref{idea1}) becomes that
\begin{eqnarray}\label{largerpart}
\mathbb{P}\Big[T(y)>e^{\sqrt{2}y(1+\delta)}\Big]& \leq & e^{-\delta_1y/2}+e^{-\delta_1y}+\exp(-c^\prime e^{\delta_1 y/2})\\
& \leq & c_6 e^{-\delta_1y/2}.\nonumber
\end{eqnarray}

By the Borel-Cantelli lemma again, we conclude that almost surely
$$\limsup_{y\rightarrow\infty}\frac{\log T(y)}{y}\leq \sqrt{2}, $$
which completes the proof of the proposition. $\square$

For $\alpha>0$ and $\beta>0$, set
$$p(z,\alpha,\beta):=\mathbb{P}\Big[\exists r\leq e^{\alpha z}: R(r)\leq m(r)-\beta z\Big].$$
\begin{proposition}\label{smalldeviation}
There exists a positive constant $C_1$, independent of $(\alpha,\beta,z)$, such that for any $z\geq z(\alpha,\beta)$,
\begin{equation}
p(z,\alpha,\beta)\leq C_1\exp\bigg(-\frac{\beta z}{6\sqrt{2}}\bigg).
\end{equation}
\end{proposition}

\noindent\textbf{Proof:} It follows from (\ref{liminf}) that as $z\rightarrow\infty$,
\begin{eqnarray*}
p(z,\alpha,\beta)& = &\mathbb{P}\bigg[\exists r\leq e^{\alpha z}: R(r)\leq m(r)-\beta z\bigg]\longrightarrow 0.
\end{eqnarray*}
Hence, there exists $z_0(\alpha,\beta)$ large enough, such that for all $z\geq z_0(\alpha,\beta)$,
\begin{equation}\label{firstestimate}
\mathbb{P}\Big[\exists r\leq e^{\alpha z}: R(r)\leq m(r)-\beta z/2\Big]\leq 1/2.
\end{equation}
\noindent For any $b_z<e^{\alpha z}$, we have
\begin{multline}\label{division}
p(z,\alpha,\beta)\leq \mathbb{P}\bigg[\exists u\in\mathcal{N}(b_z), s.t.\ \min_{s\leq b_z}X_u(s)\leq \sqrt{2}b_z-\beta z/2\bigg]\\
       + \mathbb{P}\bigg[\Big\{\exists r\leq e^{\alpha z}: R(r)\leq m(r)-\beta z\Big\}\cap\Big\{L(b_z)\geq \sqrt{2}b_z-\beta z/2\Big\}\bigg ].
\end{multline}
On the one hand, by the many-to-one lemma,
\begin{eqnarray*}
\mathbb{P}\bigg[\exists u\in\mathcal{N}(b_z), s.t.\ \min_{s\leq b_z}X_u(s)\leq \sqrt{2}b_z-\beta z/2\bigg]&\leq & \mathbb{E}\Bigg[\sum_{u\in\mathcal{N}(b_z)}\mathbf{1}_{(\min_{s\leq b_z}X_u(s)\leq \sqrt{2}b_z-\beta z/2)}\Bigg]\\
&=& e^{b_z}\mathbb{P}\bigg[ \min_{s\leq b_z} B_s \leq \sqrt{2}b_z-\beta z/2\bigg].
\end{eqnarray*}
On the other hand, by simple observations,
\begin{eqnarray*}
&&\mathbb{P}\bigg[\Big\{\exists r\leq e^{\alpha z}: R(r)\leq m(r)-\beta z\Big\}\cap\Big\{L(b_z)\geq \sqrt{2}b_z-\beta z/2\Big\}\bigg ]\\
&\leq&\mathbb{P}\bigg[\bigcap_{u\in\mathcal{N}(b_z)}\Big\{\exists t \leq e^{\alpha z}, s.t. \ R^u(t)<m(t)-\beta z/2\Big\}\bigg]\\
&=& \mathbb{E}\Bigg[\prod_{u\in\mathcal{N}(b_z)}\mathbb{P}\Big[\exists t\leq e^{\alpha z}, s.t.\ R(t)<m(t)-\beta z/2\Big]\Bigg],
\end{eqnarray*}
where the last equality follows from the branching property. Going back to (\ref{division}), one has
\begin{equation*}
p(z,\alpha,\beta)\leq e^{b_z}\mathbb{P}\bigg[ \min_{s\leq b_z} B_s \leq \sqrt{2}b_z-\beta z/2\bigg]+\mathbb{E}\Bigg[\prod_{u\in\mathcal{N}(b_z)}\mathbb{P}\Big[\exists t\leq e^{\alpha z}, s.t.\ R(t)<m(t)-\beta z/2\Big]\Bigg].
\end{equation*}
Let $b_z=\frac{\beta}{6\sqrt{2}}z$. Then, by (\ref{firstestimate}), for all $z\geq z(\alpha,\beta):=\max\{z_0(\alpha,\beta), \frac{1}{\beta}\}$,
\begin{eqnarray*}
p(z,\alpha,\beta)&\leq & e^{b_z}\mathbb{P}\Big[\min_{s\leq b_z}B_s\leq -\beta z/3\Big]+\mathbb{E}\bigg[\Big(\frac{1}{2}\Big)^{N(b_z)}\bigg]\\
       & \leq &c_7e^{-3b_z}+ e^{-b_z}\leq C_1\exp\bigg(-\frac{\beta z}{6\sqrt{2}}\bigg),
       \end{eqnarray*}
with $C_1:=c_7+1$, which completes the proof of the proposition.  $\square$

\begin{corollary}\label{cor}
For any $\delta\in(0,1)$, there exists some $s(\delta)\geq 1$, such that for all $s\geq s(\delta)$,
\begin{equation}
\mathbb{P}\Big[R(s)\leq \sqrt{2}(1-\delta)s\Big]\leq C_1\exp\bigg(-\frac{\delta s}{12\sqrt{2}}\bigg).
\end{equation}
\end{corollary}
\textbf{Proof: } Since we always have $m(s)-\delta s/2\geq \sqrt{2}(1-\delta)s$ when $s$ is sufficiently large,
$$\mathbb{P}\bigg[R(s)\leq \sqrt{2}(1-\delta)s\bigg]\leq \mathbb{P}\bigg[\exists r\leq e^s: R(r)\leq m(r)-\delta s/2\bigg].$$
which by Proposition \ref{smalldeviation} is bounded by $C_1\exp\big(-\frac{\delta s}{12\sqrt{2}}\big)$ for all $s$ large enough.  $\square$

\section{The case of two independent branching Brownian motions}
 We consider two independent branching Brownian motions, denoted by $\mathbb{X}^A(\cdot)$ and $\mathbb{X}^B(\cdot)$. Suppose that $\mathbb{P}[X^A(0)=0]=\mathbb{P}[X^B(0)=z]=1$ with $z>0$, where $\mathbb{X}^A(0)$ and $\mathbb{X}^B(0)$ represent the position of the roots, respectively. We write $R^A(\cdot)$  ($ R^B(\cdot)$, respectively) for the position of rightmost particle of the BBM $\mathbb{X}^A(\cdot)$ ($\mathbb{X}^B(\cdot)$, respectively). We define, for any $y>0$,
 \begin{eqnarray*}
 T^A(y)&:=&\inf\{t\geq 1; R^A(t)> m(t)+y\};\\
 T^B(y) &:=&\inf\{t\geq 1; R^B(t)> m(t)+y\}.
 \end{eqnarray*}
 Let $T^{A>B}$ be the first time when the rightmost point of $\mathbb{X}^A$ exceeds that of $\mathbb{X}^B$, i.e.,
 $$T^{A>B}=T^{A>B}(z):=\inf\{t\geq 0; R^A(t)>R^B(t)\}.$$
 We immediately observe that the distribution of $T^{A>B}(z)$ merely depends on the parameter $z$. Actually, we can take another pair of independent standard BBM's (both rooted at the origin), namely, $\mathbb{X}^{\Rmnum{1}}(\cdot)$ and $\mathbb{X}^{\Rmnum{2}}(\cdot)$. Their rightmost positions are denoted by $R^{\Rmnum{1}}(\cdot)$ and $R^{\Rmnum{2}}(\cdot)$, respectively. For any positive $z$, let
 $$\mathcal{T}(z):=\inf\{t\geq 0: R^{\Rmnum{1}}(t)-R^{\Rmnum{2}}(t)>z\}.$$
 Then $T^{A>B}(z)$ is distributed as $\mathcal{T}(z)$. Besides, $z\mapsto\mathcal{T}(z)$ is increasing.
\begin{proposition}\label{twobbm}
The following convergence holds almost surely
\begin{equation}\label{firstconclusion }
\lim_{z\rightarrow\infty}\frac{\log \mathcal{T}(z)}{z}=\sqrt{2}.
\end{equation}
\end{proposition}

\noindent\textbf{Proof:} For any $\delta\in(0,1)$,
\begin{equation*}
\mathbb{P}\bigg[\mathcal{T}(z)\leq e^{\sqrt{2}z(1-\delta)}\bigg]=\mathbb{P}\bigg[T^{A>B}(z)\leq e^{\sqrt{2}z(1-\delta)}\bigg]\leq p_1+p_2,
\end{equation*}
where
\begin{eqnarray*}
p_1&:=&\mathbb{P}\bigg[\exists t\leq e^{\sqrt{2}z(1-\delta)}, s.t.\  R^B(t)<m(t)+z-\delta z/2 \bigg],\\
p_2&:=&\mathbb{P}\bigg[\bigg\{T^{A>B}\leq e^{\sqrt{2}z(1-\delta)}\bigg\}\cap\bigg\{ R^B(t)\geq m(t)+(1-\delta/2)z, \forall t\leq e^{\sqrt{2}z(1-\delta)}\bigg\}\bigg].
\end{eqnarray*}
Clearly, $p_1=\mathbb{P}\bigg[\exists t\leq e^{\sqrt{2}z(1-\delta)}, s.t.\  R(t)<m(t)-\delta z/2\bigg]=p(z, \sqrt{2}(1-\delta),\delta/2)$. By Proposition \ref{smalldeviation}, for all $z\geq z(\delta)$,
\begin{equation*}
p_1\leq C_1\exp\bigg(-\frac{\delta z}{12\sqrt{2}}\bigg).
\end{equation*}
At the same time, we notice that  
\begin{multline}
\bigg\{T^{A>B}\leq e^{\sqrt{2}z(1-\delta)}\bigg\}\cap\bigg\{ R^B(t)\geq m(t)+(1-\delta/2)z, \forall t\leq e^{\sqrt{2}z(1-\delta)}\bigg\}\subset\\
\bigg\{\exists t\leq e^{\sqrt{2}z(1-\delta)}: R^A(t)\geq R^B(t)\geq m(t)+(1-\delta/2)z\bigg\}\subset\bigg\{T^A\Big((1-\delta/2)z\Big)\leq e^{\sqrt{2}z(1-\delta)}\bigg\}.
\end{multline}
\noindent This yields that
$$p_2\leq \mathbb{P}\bigg[T\Big((1-\delta/2)z\Big)\leq e^{\sqrt{2}z(1-\delta)}\bigg]\leq c_5z^2 e^{-\delta z/\sqrt{2}},$$
\noindent because of the inequality (\ref{lesspart}).

As a result,
\begin{equation}\label{twobbmlesspart}
\mathbb{P}\bigg[\mathcal{T}(z)\leq e^{\sqrt{2}z(1-\delta)}\bigg]\leq C_1\exp\bigg(-\frac{\delta z}{12\sqrt{2}}\bigg)+c_5z^2 e^{-\delta z/\sqrt{2}}\leq c_8 \exp\bigg(-\frac{\delta z}{12\sqrt{2}}\bigg),
\end{equation}
\noindent for some constant $c_8>0$ and all $z$ large enough. Thus, by the Borel-Cantelli lemma,
$$\liminf_{z\rightarrow\infty}\frac{\log \mathcal{T}(z)}{z} \geq \sqrt{2} \text{   almost surely. }$$

To prove the upper bound, we observe that
\begin{equation}\label{tbbmlargepart}
\mathbb{P}\bigg[\mathcal{T}(z)>e^{\sqrt{2}z(1+\delta)}\bigg]=\mathbb{P}\bigg[T^{A>B}(z)>e^{\sqrt{2}z(1+\delta)}\bigg]\leq q_1+q_2,
\end{equation}
where
\begin{eqnarray*}
q_1&:=& \mathbb{P}\bigg[\bigg\{T^A\Big(z(1+\delta/2)\Big)>e^{\sqrt{2}z(1+\delta)}\bigg\} \cup\bigg\{ T^A\Big(z(1+\delta/2)\Big)<e^{\sqrt{2}z}\bigg\}\bigg],\\
q_2&:=& \mathbb{P}\bigg[e^{\sqrt{2}z}\leq T^A\Big(z(1+\delta/2)\Big)\leq e^{\sqrt{2}z(1+\delta)}< T^{A>B}(z)\bigg].
\end{eqnarray*}
Notice that $T^A(y)$ is distributed as $T(y)$ for any $y>0$. According to the inequalities (\ref{lesspart}) and (\ref{largerpart}), there exists $\delta_2:=\delta_2(\delta)>0$ such that $q_1\leq e^{-\delta_2 z}$ for $z$ large enough. It remains to estimate $q_2$:
 \begin{eqnarray*}
q_2&\leq & \int_{e^{\sqrt{2}z}}^{e^{\sqrt{2}z(1+\delta)}} \mathbb{P}\bigg[T^A\Big(z(1+\delta/2)\Big)\in dr\bigg]\mathbb{P}\bigg[T^{A>B}>r\bigg\vert T^A\Big(z(1+\delta/2)\Big)=r\bigg]\\
&\leq & \int_{e^{\sqrt{2}z}}^{e^{\sqrt{2}z(1+\delta)}} \mathbb{P}\bigg[T^A\Big(z(1+\delta/2)\Big)\in dr\bigg]\mathbb{P}\bigg[R^B(r)>m(r)+z(1+\delta/2)\bigg].
\end{eqnarray*}
By the inequality (\ref{ineq}) again, this tells that
\begin{eqnarray*}
q_2&\leq &\int_{e^{\sqrt{2}z}}^{e^{\sqrt{2}z(1+\delta)}} \mathbb{P}\bigg[T^A\Big(z(1+\delta/2)\Big)\in dr\bigg]c_2 (z+1)^2 e^{-\sqrt{2}\delta z/2}\\
&\leq & c_2 (z+1)^2 e^{-\sqrt{2}\delta z/2}.
\end{eqnarray*}
Thus, recalling (\ref{tbbmlargepart}),  we obtain that for all $z$ large enough,
$$\mathbb{P}\bigg[\mathcal{T}(z)>e^{\sqrt{2}z(1+\delta)}\bigg]\leq e^{-\delta_2 z}+c_2 (z+1)^2 e^{-\sqrt{2}\delta z/2}.$$
\noindent It follows that almost surely $ \limsup_{z\rightarrow\infty}\frac{\log \mathcal{T}(z)}{z}\leq \sqrt{2}$. Proposition \ref{twobbm} is proved. $\square$

\section{Proof of Theorem \ref{leftmost particle}}

For any $k\in\mathbb{N}_+$ and $\delta\in(0, 1/20)$, we define
$$\mathcal{N}_\delta(k):=\{u\in\mathcal{N}(k): X_u(k)\leq -\sqrt{2}(1-\delta/2)k\}.$$
In order to study the asymptotic behavior of $\tau_{\ell(s)}$ for $s\in\r_+$, we first look for a lower bound for $\min_{u\in\mathcal{N}_\delta(k)}\tau_u$ and an upper bound for $\max_{u\in\mathcal{N}_\delta(k)}\tau_u$.

Recall the definitions (\ref{shiftedsubtree}) and (\ref{shiftedpositions}) of the shifted subtrees. For any particle $u\in\mathcal{N}_\delta(k)$, we use $\mathbb{X}^{u}(\cdot)$ to represent the branching Brownian motion generated by $u$ started from the time $k$. Meanwhile, we use $\mathbb{X}^{r}(\cdot)$ to represent the branching Brownian motion generated by the rightmost point at time $k$. Accordingly, the random variable $T^{u>r}$ is defined to be the first time when $u$ has a descendant exceeding all descendants of the rightmost particle at time $k$.

Considering that $T^{u>r}\leq \tau_u$ for each $u\in\mathcal{N}_\delta(k)$, one sees that
$$\mathbb{P}\bigg[\bigcup_{u\in\mathcal{N}_\delta(k)}\{\tau_u\leq e^{4k(1-10\delta)}\}\bigg]\leq p_1^\prime+p_2^\prime,$$
\noindent where
\begin{eqnarray*}
p_1^\prime&:=& \mathbb{P}\Big[R(k)\leq \sqrt{2}(1-\delta/2)k\Big],\\
p_2^\prime&:=&\mathbb{P}\Big[1_{\big( R(k)\geq \sqrt{2}(1-\delta/2)k\big)}\sum_{u\in\mathcal{N}_\delta(k)}1_{\big(T^{u>r}\leq e^{4k(1-10\delta)}\big)}\Big].
\end{eqnarray*}
Given $\mathcal{F}_k$, the BBM's $\mathbb{X}^{u}$ and $\mathbb{X}^{r}$ are independent. Then,
\begin{eqnarray}
p_2^\prime&\leq &\mathbb{E}\bigg[\sum_{u\in\mathcal{N}_\delta(k)}1_{\big(R(k)\geq \sqrt{2}(1-\delta/2)k\big)}\mathbb{P}\Big[T^{u>r}\leq e^{4k(1-10\delta)}\Big\vert\mathcal{F}_k\Big]\bigg]\nonumber\\
&= & \mathbb{E}\bigg[\sum_{u\in\mathcal{N}_\delta(k)}1_{\big(R(k)\geq \sqrt{2}(1-\delta/2)k\big)}\mathbb{P}\Big[\mathcal{T}(R(k)-X_u(k))\leq e^{4k(1-10\delta)}\Big\vert\mathcal{F}_k\Big]\bigg].\nonumber
\end{eqnarray}
By the monotonicity of $\mathcal{T}(\cdot)$, this gives that
\begin{eqnarray}
p_2^\prime & \leq & \mathbb{E}\bigg[\sum_{u\in\mathcal{N}_\delta(k)}\mathbb{P}\Big[\mathcal{T}(2\sqrt{2}k(1-\delta/2))\leq e^{4k(1-10\delta)}\Big]\bigg]\nonumber\\
&=& \mathbb{E}\bigg[\sum_{u\in\mathcal{N}_\delta(k)}1\bigg]\mathbb{P}\Big[\mathcal{T}(2\sqrt{2}k(1-\delta/2))\leq e^{4k(1-10\delta)}\Big].\nonumber
\end{eqnarray}
Using the inequality (\ref{twobbmlesspart}), for all $k$ sufficiently large,
$$p_2^\prime\leq \mathbb{E}\bigg[\sum_{u\in\mathcal{N}_\delta(k)}1\bigg] c_8 \exp\Big(-\frac{3\delta k}{2}(1-\delta/2)\Big).$$
\noindent Then by the many-to-one lemma and by (\ref{evaluationofBM}), we obtain that
\begin{equation}\label{firstsecondpart}
\begin{array}{cl}
p_2^\prime &\leq e^k\mathbb{P}[B_k\leq -\sqrt{2}k(1-\delta/2)]c_8 \exp\Big(-\frac{3\delta k}{2}(1-\delta/2)\Big)\\
&\leq  e^{-c_9 \delta k},
\end{array}
\end{equation}
where $c_9$ is a positive constant independent of $(\delta, k)$.

In view of Corollary \ref{cor}, for large $k$, one has
\begin{equation}\label{firstfirstpart}
p_1^\prime\leq C_1\exp\bigg(-\frac{\delta k}{24\sqrt{2}}\bigg).
\end{equation}
Combining (\ref{firstfirstpart}) with (\ref{firstsecondpart}) yields that for $k$ large enough,
\begin{equation*}
\mathbb{P}\bigg[\bigcup_{u\in\mathcal{N}_\delta(k)}\{\tau_u\leq e^{4k(1-10\delta)}\}\bigg]\leq C_1\exp\bigg(-\frac{\delta k}{24\sqrt{2}}\bigg)+e^{-c_9 \delta k}.
\end{equation*}
By the Borel-Cantelli lemma, almost surely,
\begin{equation}\label{fristpart}
\liminf_{k\rightarrow\infty}\frac{\log \min_{u\in\mathcal{N}_\delta(k)}\tau_u}{k}\geq 4(1-10\delta),
\end{equation}
which gives the lower bound for $\min_{u\in\mathcal{N}_\delta(k)}\tau_u$.

To obtain an upper bound for $\max_{u\in\mathcal{N}_\delta(k)}\tau_u$, let us estimate $\mathbb{P}\Big[\cup_{u\in\mathcal{N}_\delta(k)}\Big\{\tau_u\geq e^{4k(1+10\delta)}\Big\}\Big]$. We consider the subtree generated by any particle $u\in\mathcal{N}_\delta(k)$. Recall that the shifted positions of its descendants are denoted by
$$X^u_v(\cdot):=X_v(\cdot+k)-X_{u}(k)\text{ for any }v\in\mathcal{N}(\cdot+k)\text{ satisfying }u<v,$$
and that $R^u(\cdot):=\max X^u_v(\cdot)$. We set $T^u(y):=\inf\{t\geq1; R^u(t)-m(t)>y\}$ for any $y>0$, which is obviously distributed as $T(y)$. Let $y=2\sqrt{2}k(1+\delta/2)$, then
\begin{equation}\label{secondpart}
\mathbb{P}\bigg[\bigcup_{u\in\mathcal{N}_\delta(k)}\Big\{\tau_u\geq e^{4k(1+10\delta)}\Big\}\bigg]\leq q_1^\prime+q_2^\prime+q_3^\prime,
\end{equation}
where
\begin{eqnarray*}
q_1^\prime&:= & \mathbb{P}\bigg[\bigcup_{u\in\mathcal{N}_\delta(k)}\bigg(\Big\{T^u(y)\geq e^{4k(1+10\delta)}\Big\}\cup\Big\{T^u(y)\leq e^{k}\Big\}\bigg)\bigg],\\
q_2^\prime&:= &\mathbb{P}\bigg[L(k)\leq -\sqrt{2}k\bigg],\\
q_3^\prime&:= &\mathbb{P}\bigg[\bigcup_{u\in\mathcal{N}_\delta(k)}\{e^{k}<T^u(y)<e^{4k(1+10\delta)}\leq \tau_u\}; L(k)>-\sqrt{2}k\bigg].
\end{eqnarray*}
First, we observe that
$$q_1^\prime\leq \mathbb{E}\bigg[\sum_{u\in\mathcal{N}_\delta(k)}1\bigg]\mathbb{P}\bigg[\Big\{T(y)\geq e^{4k(1+10\delta)}\Big\}\cup\Big\{T(y)\leq e^{k}\Big\}\bigg].$$
Using the many-to-one lemma for the first term on the right-hand side,
$$q_1^\prime\leq e^k\mathbb{P}\bigg[B_k\leq-\sqrt{2}k(1-\delta/2)\bigg]\mathbb{P}\bigg[\Big\{T(y)\geq e^{4k(1+10\delta)}\Big\}\cup\Big\{T(y)\leq e^{k}\Big\}\bigg].$$
According to the inequalities (\ref{evaluationofBM}) (\ref{lesspart}) and (\ref{largerpart}), there exists $\delta_4:=\delta_4(\delta)>0$ such that $q_1^\prime\leq e^{-\delta_4 k}$ for $k$ large enough. Meanwhile, by (\ref{ineq}), $q_2^\prime\leq 2c_2(\log k+1)^2 k^{-3/2}$.

It remains to bound $q_3^\prime$. Since $T^u(y)$ is independent of $\mathcal{F}_k$, it follows that
\begin{eqnarray*}
q_3^\prime&\leq&\mathbb{E}\left[\sum_{u\in\mathcal{N}_\delta(k)}\int_{e^k}^{e^{4k(1+10\delta)}}\mathbb{P}\bigg[T^u(y)\in dr\bigg]\mathbb{P}\bigg[\tau_u>r; L(k)\geq-\sqrt{2}k\bigg\vert T^u(y)=r, \mathcal{F}_k\bigg]\right].
\end{eqnarray*}
Given $\{T^u(y)=r\}$ and $\mathcal{F}_k$, the event $\{\tau_u>r\}\cap\{L(k)\geq-\sqrt{2}k\}$ implies that $\cup_{w\in\mathcal{N}(k)\setminus\{u\}}\{R^\omega(r)+X_w(k)>R^u(r)+X_u(k)\geq m(r)+y-\sqrt{2}k\}$, whose probability is less than $\sum_{w\in\mathcal{N}(k)\setminus\{u\}}c_1\Big(1+(y-\sqrt{2}k-X_w(k))^2_+\Big)e^{-\sqrt{2}y+2k+\sqrt{2}X_w(k)}$ (see (\ref{ineq})). This yields that
\begin{eqnarray*}
q_3^\prime&\leq& \mathbb{E}\left[\sum_{u\in\mathcal{N}_\delta(k)}\int_{e^{k}}^{e^{4k(1+10\delta)}}\mathbb{P}\Big[T^u(y)\in dr\Big]\sum_{w\in\mathcal{N}(k)\setminus\{u\}}c_2(y+1)^2e^{-\sqrt{2}y+2k+\sqrt{2}X_w(k)}\right]\\
&\leq& \mathbb{E}\left[\sum_{u\in\mathcal{N}_\delta(k)}\sum_{w\in\mathcal{N}(k)\setminus\{u\}}c_2(y+1)^2e^{-\sqrt{2}y+2k+\sqrt{2}X_\omega(k)}\right]\\
&= &c_2(y+1)^2e^{-\sqrt{2}y+2k} \mathbb{E}\bigg[\sum_{u\in\mathcal{N}_\delta(k)}\sum_{w\in\mathcal{N}(k)\setminus\{u\}}e^{\sqrt{2}X_w(k)}\bigg].
\end{eqnarray*}

By integrating with respect to the last time at which the most recent common ancestor of $u$ and $\omega$ was alive (see e.g. \cite{Harris2011} for more details), $\mathbb{E}\Big[\sum_{u\in\mathcal{N}_\delta(k)}\sum_{\omega\neq u}e^{\sqrt{2}X_\omega(k)}\Big]$ is equal to
\begin{eqnarray*}
&&2\int_{0}^k e^{2k-s}ds\int_{\mathbb{R}}\mathbb{P}\bigg[B_s\in dx\bigg]\mathbb{P}\bigg[B_k\leq -\sqrt{2}(1-\delta/2)k\Big\vert B_s=x\bigg]\mathbb{E}\bigg[e^{\sqrt{2}B_k}\Big\vert B_s=x\bigg]\\
&=&2\int_0^k e^{2k-s}ds\int_{\mathbb{R}}\mathbb{P}\bigg[B_s\in dx\bigg]\mathbb{P}\bigg[B_k\leq -\sqrt{2}(1-\delta/2)k\Big\vert B_s=x\bigg] e^{\sqrt{2}x}e^{k-s},
\end{eqnarray*}
where the second equivalence follows from the Markov property of Brownian Motion. We rearrange the integration as follows:
\begin{eqnarray*}
\mathbb{E}\bigg[\sum_{u\in\mathcal{N}_\delta(k)}\sum_{w\in\mathcal{N}(k)\setminus\{u\}}e^{\sqrt{2}X_w(k)}\bigg]&=&2\int_0^k e^{3k-2s}\mathbb{E}\bigg[e^{\sqrt{2}B_s}; B_k\leq -\sqrt{2}(1-\delta/2)k\bigg]ds\\
&=&2\int_0^k e^{3k-2s}\mathbb{E}\bigg[e^{\sqrt{2}B_k}e^{-\sqrt{2}(B_k-B_s)}; B_k\leq -\sqrt{2}(1-\delta/2)k\bigg]ds\\
&\leq &2\int_0^k e^{3k-2s} e^{-2(1-\delta/2)k}\mathbb{E}\bigg[e^{-\sqrt{2}(B_k-B_s)}\bigg]ds,
\end{eqnarray*}
which is bounded by $e^{(2+\delta)k}$ by simple computation. Thus, $q_3^\prime\leq c_{10} k^2e^{-\delta k}$ for some constant $c_{10}>0$.

Going back to (\ref{secondpart}),
$$\mathbb{P}\bigg[\bigcup_{u\in\mathcal{N}_\delta(k)}\Big\{\tau_u\geq e^{4k(1+10\delta)}\Big\}\bigg]\leq e^{-\delta_4 k}+2c_2(\log k+1)^2 k^{-3/2}+ c_{10}k^2 e^{-\delta k},$$
for all $k$ sufficiently large.

Therefore, by the Borel-Cantelli lemma,
$$\limsup_{\mathbb{N}\owns k\rightarrow\infty}\frac{\log \max_{u\in\mathcal{N}_\delta(k)}\tau_{u}}{k}\leq 4(1+10\delta) \ \ \text{ almost surely.}$$

We now turn to study $\{\tau_{\ell(s)}; s\geq 0\}$.

On the one hand, for any $\delta>0$, we claim that almost surely for $s$ large enough, the leftmost particle $\ell(s)$ at time $s$ must have at least one descendant belonging to $\mathcal{N}_\delta(\lfloor s\rfloor+1)$.

In fact, let us write $\Upsilon_k:=\{\exists u\in\mathcal{N}(k+1): u\notin\mathcal{N}_\delta(k+1); \exists s\in[k,k+1], X_u(s)\leq -\sqrt{2}s+\delta^\prime s\}$ with $\delta^\prime:=(\sqrt{2}-1)\delta/2$. By the many-to-one lemma, we get that for $k\geq 100/\delta$,
\begin{equation*}
\mathbb{P}[\Upsilon_k]\leq \frac{4}{\delta k\sqrt{2\pi}}e^{1+k-\frac{\delta^2k^2}{8}},
\end{equation*}
which is summable over $k$. It follows that
\begin{equation}\label{upsilonestimation}
\p[\Upsilon_k \text{ infinitely often}]=0.
\end{equation}
In view of (\ref{limsup}), when $s$ is large enough, $L(s)$ always lie below $-\sqrt{2}s+\delta^\prime s$ almost surely. Combining with (\ref{upsilonestimation}), we obtain that almost surely for $k$ sufficiently large,
\begin{equation*}
\max_{s\in[k,k+1]}\tau_{\ell(s)}\leq \max_{u\in\mathcal{N}_\delta(k+1)}\tau_u+1.
\end{equation*}

On the other hand, using similar arguments, one can say that almost surely for $s$ sufficiently large, the leftmost located particle $\ell(s)$ at time $s$ must come from one particle in $\mathcal{N}_\delta(\lfloor s\rfloor)$. This gives that almost surely for $k$ sufficiently large,
\begin{equation*}
\min_{s\in[k,k+1]}\tau_{\ell(s)}\geq \min_{u\in\mathcal{N}_\delta(k)}\tau_u-1.
\end{equation*}
Thus we conclude that almost surely
$$\lim_{s\rightarrow\infty}\frac{\log \tau_{\ell(s)}}{s}=4.\qquad \square$$

\section{Proof of Theorem \ref{mainconclusion}}
It suffices to show that almost surely $\lim_{\mathbb{N}\owns k\rightarrow\infty}\frac{\log \Theta_k}{k}= 2+2\sqrt{2}$, as the sequence $\{\Theta_s; s>0\}$ is monotone.

\subsection{The lower bound of Theorem \ref{mainconclusion}} This subsection is devoted to checking that: almost surely,

$$\liminf_{k\rightarrow\infty}\frac{\log \Theta_k}{k}\geq 2+2\sqrt{2}.$$

For $0<a<\sqrt{2}$, we define
\begin{equation*}
\mathcal{Z}_a(k):=\big\{u\in\mathcal{N}(k); X_u(k)\leq -a k\big\}\text{ and }Z_a(k):=\# \mathcal{Z}_a(k).
\end{equation*}

For $0<\varepsilon<(1-\frac{a^2}{2})/2$ and $0<\delta<1$, we denote
\begin{eqnarray*}
E_k&:=&\Big\{Z_a(k)\geq \exp[k(1-\frac{a^2}{2}-\varepsilon)]\Big\},\\
D_k&:=&\Big\{\Theta_k\leq \exp[(2+2\sqrt{2}-\delta)k]\Big\}.
\end{eqnarray*}

Let us estimate $\mathbb{P}[D_k\cap E_k]$.

For any $s>0$ and $\beta>0$, we write $\Gamma=\Gamma(s,\beta):=\big\{N(s)=1, L(s)\leq -\beta s\big\}$. Similarly, let $\Gamma_u:=\big\{N^u(s)=1, L^u(s)\leq -\beta s\big\}$ for every $u\in\mathcal{N}(k)$. Then,
\begin{equation}
\mathbb{P}\bigg[D_k\cap E_k\bigg]\leq \mathbb{P}\bigg[\bigg(\bigcap_{u\in\mathcal{Z}_a(k)}\Gamma_u^c\bigg)\cap E_k\bigg]+\mathbb{P}\bigg[\bigg(\bigcup_{u\in\mathcal{Z}_a(k)}\Gamma_u\bigg)\cap D_k\bigg]\label{lowerbound}.
\end{equation}
By the branching structure, we obtain that
\begin{equation}
\mathbb{P}\bigg[\bigg(\bigcap_{u\in\mathcal{Z}_a(k)}\Gamma_u^c\bigg)\cap E_k\bigg]\leq \mathbb{P}\bigg[\bigg(1-\mathbb{P}[\Gamma]\bigg)^{Z_a(k)}; E_k\bigg]\leq e^{-\mathbb{P}[\Gamma]\exp[k(1-\frac{a^2}{2}-\varepsilon)]}.
\end{equation}
Clearly, $\mathbb{P}[\Gamma]=e^{-s}\mathbb{P}[B_s\leq -\beta s]$. By (\ref{convergenceofBM}), one sees that, for $\varepsilon>0$ small and $s$ large enough,
\begin{equation}\label{lowerbound1}
\mathbb{P}\bigg[\bigg(\bigcap_{u\in\mathcal{Z}_a(k)}\Gamma_u^c\bigg)\cap E_k\bigg]\leq \exp\bigg\{-\exp\Big[-s(1+\frac{\beta^2}{2}+\varepsilon)+k(1-\frac{a^2}{2}-\varepsilon)\Big]\bigg\},
\end{equation}
which is bounded by $e^{-e^{k\varepsilon}}$ if we choose $s=\frac{1-\frac{a^2}{2}-2\varepsilon}{1+\frac{\beta^2}{2}+\varepsilon}k$ with $k$ sufficiently large.

It remains to bound $\Omega:=\mathbb{P}\bigg[\Big(\bigcup_{u\in\mathcal{Z}_a(k)}\Gamma_u\Big)\cap D_k\bigg]$ for $s=\frac{1-\frac{a^2}{2}-2\varepsilon}{1+\frac{\beta^2}{2}+\varepsilon}k$. Recalling the definition of $\Theta_k$, one sees that for any $\rho\in(0,2)$,
\begin{multline}\label{omegabound}
\Omega\leq \mathbb{P}\bigg[\bigcup_{u\in\mathcal{Z}_a(k)}\Big\{\tau_u<e^{\rho k}\Big\}\bigg]\\
+\mathbb{P}\bigg[\bigg(\bigcup_{u\in\mathcal{Z}_a(k)}\Gamma_u\bigg)\cap\bigg(\bigcap_{u\in\mathcal{Z}_a(k)}\Big\{e^{\rho k}\leq\tau_u\leq e^{(2+2\sqrt{2}-\delta) k}\Big\}\bigg)\bigg]=:\Omega_a+\Omega_b.
\end{multline}

We choose now $\rho=1-2\varepsilon$ and $z=(\sqrt{2}-\frac{a^2}{2\sqrt{2}}-\frac{\varepsilon}{\sqrt{2}})k$. Then comparing $T^u(z)$ and $e^{\rho k}$ for every $u\in\mathcal{Z}_a(k)$ tells us that
\begin{equation}\label{omegaa}
\Omega_a\leq \mathbb{P}\bigg[\bigcup_{u\in\mathcal{Z}_a(k)}\{T^u(z)<e^{\rho k}\}\bigg]+\mathbb{P}\bigg[\bigcup_{u\in\mathcal{Z}_a(k)}\{\tau_u <e^{\rho k}\leq T^u(z)\}\bigg].
\end{equation}
It follows from the branching property that the first term of the right-hand side is bounded by $\mathbb{E}[Z_a(k)]\mathbb{P}[T(z)<e^{\rho k}]$, which is $e^k\mathbb{P}[B_k\leq -ak]\mathbb{P}[T(z)<e^{\rho k}]$ by the many-to-one lemma. In view of the inequalities (\ref{evaluationofBM}) and (\ref{lesspart}), one immediately has
\begin{equation}\label{omegaa1}
\mathbb{P}\bigg[\bigcup_{u\in\mathcal{Z}_a(k)}\{T^u(z)<e^{\rho k}\}\bigg]\leq e^k\frac{\sqrt{k}}{ak}e^{-a^2k/2}c_5 e^{\rho k}z^2 e^{-\sqrt{2}z}\leq e^{-\eta k},
\end{equation}
for some $\eta:=\eta(\varepsilon)>0$ small enough.

For the second term of the right-hand side in (\ref{omegaa}), we observe that for any $u\in\mathcal{Z}_a(k)$, $\{\tau_u <e^{\rho k}\leq T^u(z)\}$ implies that at time $\tau_u<e^{\rho k}$, the rightmost position $R(k+\tau_u)$ is exactly equal to $R^u(\tau_u)+X_u(k)$, which is less than $m(\tau_u)+z-ak$. Hence, the event $\cup_{u\in\mathcal{Z}_a(k)}\{\tau_u <e^{\rho k}\leq T^u(z)\}$ ensures that there exists some time $r<e^{\rho k}$ such that the rightmost position $R(k+r)$ is less than $m(r)+z-ak$. This gives that
\begin{equation}
\mathbb{P}\bigg[\bigcup_{u\in\mathcal{Z}_a(k)}\{\tau_u <e^{\rho k}\leq T^u(z)\}\bigg]\leq \mathbb{P}\bigg[\exists r\leq e^{\rho k}\text{ s.t. }R(k+r)\leq m(r)+z-ak\bigg].
\end{equation}
Notice that with our choice of $\rho$ and $z$, Proposition \ref{smalldeviation} can be applied to show that for all $k$ sufficiently large,
\begin{equation}
\mathbb{P}\bigg[\bigcup_{u\in\mathcal{Z}_a(k)}\{\tau_u <e^{\rho k}\leq T^u(z)\}\bigg]\leq e^{-\eta k}.
\end{equation}
Combined with (\ref{omegaa1}), the inequality (\ref{omegaa}) becomes $\Omega_a\leq 2e^{-\eta k}$.

As shown in (\ref{omegabound}), it remains to study $\Omega_b$. For the particles $u\in\mathcal{Z}_a(k)$ such that $N^u(s)=1$, we focus on the subtree rooted at $u$ but started from time $k+s$. Define
$$\widetilde{R}^u(t):=\max\Big\{X_v(k+s+t)-X_u(k+s); v\in\mathcal{N}(k+s+t), u<v\Big\},\forall t\geq0;$$
and
$$\widetilde{T}^u(y):=\inf\Big\{t\geq 1; \widetilde{R}^u(t)\geq m(t)+y\Big\},\forall y>1.$$
Since $(\widetilde{R}^u(t),t\geq0)$ is distributed as $(R(t),t\geq0)$, $\widetilde{T}^u(y)$ has the same law as $T(y)$. Let us take $\sqrt{2}x=k(2+2\sqrt{2}-\delta/2)$. Comparing $\widetilde{T}^u(x)$ with $e^{k(2+2\sqrt{2}-\delta)}$ yields that
\begin{multline}
\Omega_b\leq \mathbb{P}\bigg[\exists \ \omega\in\mathcal{Z}_a(k),\text{ s.t. }N^\omega(s)=1, L^\omega(s)\leq -\beta s,\widetilde{T}^\omega(x)\leq e^{k(2+2\sqrt{2}-\delta)}\bigg] \\
+\mathbb{P}\bigg[\exists u\in\mathcal{Z}_a(k)\text{ s.t. } N^u(s)=1,L^u(s)\leq -\beta s,e^{(1-2\varepsilon)k}\leq \tau_u\leq e^{k(2+2\sqrt{2}-\delta)}<\widetilde{ T}^u(x)\bigg]=:\Omega_{b1}+\Omega_{b2}.
\end{multline}

By first conditioning on $\mathcal{F}_{k+s}$ and then on $\mathcal{F}_{k}$, one has
\begin{eqnarray}\label{omegab1}
\Omega_{b1}& \leq&  \mathbb{E}\Big[Z_a(k)\Big]\mathbb{P}\Big[\Gamma\Big]\mathbb{P}\Big[T(x)\leq e^{k(2+2\sqrt{2}-\delta)}\Big]\\
&=& e^k\mathbb{P}\Big[B_k\leq -ak\Big]e^{-s}\mathbb{P}\Big[B_s\leq -\beta s\Big]\mathbb{P}\Big[T(x)\leq e^{k(2+2\sqrt{2}-\delta)}\Big].\nonumber
\end{eqnarray}
For $0<\varepsilon<\min\{\delta/8,(1-\frac{a^2}{2})/2\}$, by (\ref{lesspart}) and (\ref{evaluationofBM}),
\begin{eqnarray}\label{omegab1bound}
\Omega_{b1}& \leq & e^{3\varepsilon k}c_5 x^2e^{-\sqrt{2}x}e^{k(2+2\sqrt{2}-\delta)}\leq c_{11} k^2e^{-\varepsilon k}.
\end{eqnarray}
On the other hand, the event $\{\exists u\in\mathcal{Z}_a(k)\text{ s.t. } N^u(s)=1,L^u(s)\leq -\beta s,e^{(1-2\varepsilon)k}\leq \tau_u\leq e^{k(2+2\sqrt{2}-\delta)}<\widetilde{ T}^u(x)\}$ implies that there exists some time $r\in[e^{(1-2\varepsilon) k}-s, e^{k(2+2\sqrt{2}-\delta)}-s]$ such that the rightmost position $R(k+s+r)$ is less than $-ak-\beta s+m(r)+x$. Thus,
\begin{eqnarray*}
\Omega_{b2}&=&\mathbb{P}\bigg[\exists u\in\mathcal{Z}_a(k)\text{ s.t. } N^u(s)=1,L^u(s)\leq -\beta s,e^{(1-2\varepsilon)k}\leq \tau_u\leq e^{k(2+2\sqrt{2}-\delta)}<\tilde{ T}^u(x)\bigg]\\
& \leq &\mathbb{P}\bigg[\exists r\in[e^{(1-2\varepsilon) k}, e^{k(2+2\sqrt{2}-\delta)}+k],\text{ s.t. }R(r)\leq m(r-k-s)+x-ak-\beta s\bigg].
\end{eqnarray*}
By taking $a=\beta=2-\sqrt{2}$, we obtain $s=\frac{1-\frac{a^2}{2}-2\varepsilon}{1+\frac{\beta^2}{2}+\varepsilon}k=(\frac{1}{\sqrt{2}}-\varepsilon_1)k$ for some sufficiently small $\varepsilon_1=\varepsilon_1(\varepsilon)>0$. Let $\delta\geq 8\sqrt{2}\varepsilon_1$, then $m(r-k-s)+x-ak-\beta s\leq m(r)-\varepsilon_1 k$. By Proposition \ref{smalldeviation}, for $k$ large enough,
\begin{equation}\label{omegab2bound}
\Omega_{b2}\leq C_1e^{-\varepsilon_1 k/6\sqrt{2}}.
\end{equation}
Since $\Omega_b\leq\Omega_{b1}+\Omega_{b2}$, it follows from (\ref{omegab1bound}) and (\ref{omegab2bound}) that $\Omega_b\leq c_{11} k^2e^{-\varepsilon k}+C_1e^{-\varepsilon_1 k/6\sqrt{2}}$. Combined with the fact that $\Omega_a\leq 2 e^{-\eta k}$,  (\ref{omegabound}) implies that
\begin{eqnarray*}
\Omega_1&\leq& \Omega_a+\Omega_{b}\leq 2e^{-\eta k}+c_{11} k^2e^{-\varepsilon k}+C_1e^{-\varepsilon_1 k/6\sqrt{2}}.
\end{eqnarray*}
\noindent According to the inequality (\ref{lowerbound}), for $0<\varepsilon<\delta/8$, $\eta(\varepsilon)>0$, $0<\varepsilon_1(\varepsilon)\leq \delta/8\sqrt{2}$ with $\delta$ sufficiently small, and for all $k$ sufficiently large,
\begin{eqnarray*}
\mathbb{P}[D_k\cap E_k]&\leq& \exp\Big[-e^{\big(k(1-\frac{a^2}{2}-\varepsilon)-s(1+\frac{\beta^2}{2}+\varepsilon)\big)}\Big]+\Omega_1\\
&\leq &e^{-e^{\varepsilon k}}+2e^{-\eta k}+c_{11} k^2e^{-\varepsilon k}+C_1e^{-\varepsilon_1 k/6\sqrt{2}}.
\end{eqnarray*}
\noindent Consequently,
$$\sum_k\mathbb{P}[D_k\cap E_k]< \infty.$$
By Borel-Cantelli lemma, $\mathbb{P}[D_k\cap E_k \text{ i.o.}]=0$. Recall that $E_k:=\Big\{Z_a(k)\geq \exp\Big[k(1-\frac{a^2}{2}-\varepsilon)\Big]\Big\}.$
Biggins \cite{Biggins1977} showed that almost surely
\begin{equation}\label{Biggins}
\lim_{k\rightarrow\infty}\frac{\log Z_a(k)}{k}=1-\frac{a^2}{2}.
\end{equation}

\noindent Therefore, for any $\delta>0$ small, $\liminf_{k\rightarrow\infty}\frac{\log \Theta_k}{k}\geq 2+2\sqrt{2}-\delta$ almost surely. This implies the lower bound of Theorem \ref{mainconclusion}. $\square$

\subsection{The upper bound of Theorem \ref{mainconclusion}}
It remains to prove the upper bound, namely, almost surely,
\begin{equation}
\limsup_{k\rightarrow\infty}\frac{\log\Theta_k}{k}\leq 2+2\sqrt{2}.
\label{upperbound of main thm}
\end{equation}

Before bringing out the proof of (\ref{upperbound of main thm}), let us state some preliminary results first.

For $M\in\mathbb{N}_+$, define
\begin{equation*}
\sigma_M:=\inf\{s>0; N(s)=M+1\}.
\end{equation*}
Clearly, $\sigma_M$ is a stopping time with respect to $\{\mathcal{F}_s;\;s\geq0\}$. Since $N(s)$ follows the geometric distribution, one sees that for any $s\geq 0$, $\mathbb{P}[\sigma_M\leq s]=\mathbb{P}[N(s)\geq 1+M]=(1-e^{-s})^{M}$. Moreover, $\sigma_M$ has a density function, denoted by $f_M$, as follows:
\begin{equation}\label{density}
f_M(s):=1_{(s\geq 0)}Me^{-s}(1-e^{-s})^{M-1}\leq 1_{(s\geq 0)}Me^{-s}.
\end{equation}
Recall that $L(s)=\inf_{u\in\mathcal{N}(s)}X_u(s)$. Let $L(\sigma_M)$ denote the leftmost position at time $\sigma_M$. Notice that at time $\sigma_M$, there are $M+1$ particles which occupy at most $M$ different positions. This tells us that for any $s$, $ \mu>0$,
\begin{equation}\label{condlaw}
\mathbb{P}\Big[L(s)\leq-\mu s\big\vert \sigma_M=s\Big]\leq M\mathbb{P}\Big[B_s\leq-\mu s\Big]\leq\frac{M}{\mu s}e^{-\mu^2 s/2},
\end{equation}
where the last inequality holds because of (\ref{evaluationofBM}).

Let $\varepsilon\in (0,1/2)$. For $r>1/\varepsilon$ and $0<s<r$, we set $\lambda:=\lambda(s,r)>0$ such that $s(1+\frac{\lambda^2}{2})=r$. Let
\begin{eqnarray}
\Phi(r,\lambda)&:=&\{\sigma_M>r-1\}\cup\{\varepsilon r\leq \sigma_M\leq r-1, L(\sigma_M)\leq-\lambda(\sigma_M, r)\sigma_M\},\\
\Psi(r,\lambda)&:=&\big\{\varepsilon r\leq \sigma_M \leq r-1, L(\sigma_M)\geq- \lambda(\sigma_M,r) \sigma_M \big\}.
\end{eqnarray}
We have the following lemma, which gives some results of the random vector $(\sigma_M,\, L(\sigma_M))$.
\begin{lemma}\label{jointdis}
\begin{description}
\item[ (i)] There exists a constant $c_{12}>0$ such that
\begin{equation}\label{phiupp}
\mathbb{P}\bigg[\Phi\Big(r,\lambda\Big)\bigg]\leq c_{12} M^2 re^{-r}.
\end{equation}
\item[(ii)]There exists a constant $c_{13}>0$ such that
\begin{equation}\label{psiupp}
\mathbb{E}\Bigg[e^{2\sigma_M-\sqrt{2}L(\sigma_M)};\Psi(r,\lambda)\Bigg]\leq c_{13}M^2 r^2e^{\sqrt{2}r}.
\end{equation}
\end{description}
\end{lemma}

\textit{Proof of Lemma \ref{jointdis}.} (i) Observe that
\begin{eqnarray*}
\mathbb{P}\bigg[\Phi\Big(r,\lambda\Big)\bigg]&\leq &\mathbb{P}\bigg[\sigma_M>r-1\bigg]+\mathbb{P}\bigg[\varepsilon r\leq \sigma_M\leq r-1, L(\sigma_M)\leq-\lambda(\sigma_M,r)\sigma_M\bigg]\\
&=&\int_{r-1}^\infty f_M(s)ds+\int_{\varepsilon r}^{r-1}\mathbb{P}\bigg[ L(s)\leq -\lambda(s,r) s\Big\vert \sigma_M=s\bigg]f_M(s)ds\\
&\leq &\int_{r-1}^\infty Me^{-s}ds+\int_{\varepsilon r}^{r-1}\frac{M}{\lambda(s,r) s}e^{-\lambda(s, r)^2 s/2} Me^{-s}ds,
\end{eqnarray*}
where the last inequality follows from (\ref{density}) and (\ref{condlaw}). A few lines of simple computation yield (\ref{phiupp}).

(ii) Let us prove the inequality (\ref{psiupp}). By Fubini's theorem, we rewrite the expectation $\mathbb{E}\Bigg[e^{2\sigma_M-\sqrt{2}L(\sigma_M)};\Psi(r,\lambda)\Bigg]$ as follows:
\begin{eqnarray}\label{rewriteintegration}
&&\int_{\varepsilon r}^{r-1}e^{2s}f_M(s)\mathbb{E}\bigg[e^{-\sqrt{2}L(s)}; L(s)\geq-\lambda(s,r)s\bigg\vert\sigma_M=s\bigg]ds\\
&=&\int_{\varepsilon r}^{r-1}e^{2s}f_M(s)\int_{-\lambda(s, r) s}^{+\infty}\sqrt{2}e^{-\sqrt{2}x}\mathbb{P}\bigg[-\lambda(s, r) s\leq L(s)\leq x\bigg\vert\sigma_M=s\bigg]dxds\nonumber.
\end{eqnarray}
For $\varepsilon r\leq s\leq r$, one sees that $\lambda(s, r)=\sqrt{2(r-s)/s}\geq \sqrt{2/(r-1)}>0$. We choose $0<\lambda_0=\min\{1-\sqrt{2}/2,\,\sqrt{2/(r-1)}\}$ so that
\begin{eqnarray}
&&\int_{-\lambda(s, r) s}^{+\infty}\sqrt{2}e^{-\sqrt{2}x}\mathbb{P}\bigg[-\lambda(s, r) s\leq L(s)\leq x\bigg\vert\sigma_M=s\bigg]dx\nonumber\\
&\leq& \int_{-\lambda_0 s}^{+\infty} \sqrt{2}e^{-\sqrt{2}x}dx+\int_{-\lambda(s, r) s}^{-\lambda_0 s}\sqrt{2}e^{-\sqrt{2}x}\mathbb{P}\bigg[-\lambda(s, r) s\leq L(s)\leq x\bigg\vert\sigma_M=s\bigg]dx\nonumber\\
&\leq& e^{(\sqrt{2}-1)s}+\int_{-\lambda(s, r) s}^{-\lambda_0 s}\sqrt{2}e^{-\sqrt{2}x}\mathbb{P}\bigg[-\lambda(s, r) s\leq L(s)\leq x\bigg\vert\sigma_M=s\bigg]dx.
\label{cutting}
\end{eqnarray}
The last term on the right-hand side of (\ref{cutting}), by a change of variable $x=-\mu s$, becomes
\begin{eqnarray*}
&&\int_{\lambda_0}^{\lambda(s, r)} \sqrt{2}s e^{\sqrt{2}\mu s}\mathbb{P}\bigg[-\lambda(s, r) s\leq L(s)\leq -\mu s\bigg\vert\sigma_M=s\bigg]d\mu\\
&\leq& \int_{\lambda_0}^{\lambda(s, r)} \sqrt{2}s e^{\sqrt{2}\mu s}\mathbb{P}\bigg[L(s)\leq -\mu s\bigg\vert\sigma_M=s\bigg]d\mu\\
&\leq& \int_{\lambda_0}^{\lambda(s, r)} \sqrt{2}s e^{\sqrt{2}\mu s} \frac{M}{\mu s}e^{-\mu^2 s/2}d\mu,
\end{eqnarray*}
where the last inequality comes from (\ref{condlaw}).
Going back to (\ref{rewriteintegration}),
\begin{eqnarray}
\mathbb{E}\bigg[e^{2\sigma_M-\sqrt{2}L(\sigma_M)};\Psi(r,\lambda)\bigg]&\leq &\int_{\varepsilon r}^{r-1}e^{2s}f_M(s)\bigg(\int_{\lambda_0}^{\lambda(s, r)} \sqrt{2s}\frac{M}{\mu}e^{(\sqrt{2}\mu-\frac{\mu^2}{2})s}d\mu+e^{(\sqrt{2}-1) s}\bigg)ds.\nonumber
\end{eqnarray}
By (\ref{density}), this is bounded by
\begin{equation*}
\int_{\varepsilon r}^{r-1}\int_{\lambda_0}^{\lambda(s, r)} \sqrt{2s}\frac{M^2}{\mu}\exp\Big\{(1+\sqrt{2}\mu-\frac{\mu^2}{2})s\Big\}d\mu ds+M \int_{\varepsilon r}^{r-1} e^{\sqrt{2}s}ds.
\end{equation*}
Notice that if $s\leq r/2$, then $\lambda(s, r)=\sqrt{2(\frac{r}{s}-1)}\geq\sqrt{2}$. It follows that
\begin{equation}
\max_{0<\mu\leq\lambda(s, r)}(1+\sqrt{2}\mu-\frac{\mu^2}{2})s= 2s\leq r<\sqrt{2}r.
\end{equation}
Otherwise, $r/2<s<r$ implies $\lambda(s, r)<\sqrt{2}$. Hence, $\max_{0<\mu\leq\lambda(s, r)}(1+\sqrt{2}\mu-\frac{\mu^2}{2})s$ is achieved when $\mu=\lambda(s,r)$, which equals
\begin{equation}
\Big(1+\sqrt{2}\lambda(s,r)-\frac{\lambda(s,r)^2}{2}\Big)s=\Big(\frac{1+\sqrt{2}\lambda(s,r)-\frac{\lambda(s,r)^2}{2}}{1+\frac{\lambda(s,r)^2}{2}}\Big)r.
\end{equation}
It is bounded by $\sqrt{2}r$ since $\max_{z\geq0}\frac{1+\sqrt{2}z-z^2/2}{1+z^2/2}=\sqrt{2}$. Combing the two cases, we obtain that
\begin{equation*}
\max_{\epsilon r\leq s\leq r-1; 0<\mu\leq\lambda(s, r)}e^{(1+\sqrt{2}\mu-\frac{\mu^2}{2})s}\leq e^{\sqrt{2}r}.
\end{equation*}
This implies that
\begin{eqnarray*}
\mathbb{E}\bigg[e^{2\sigma_M-\sqrt{2}L(\sigma_M)};\Psi(r,\lambda)\bigg]&\leq& \int_{\varepsilon r}^{r-1}\int_{\lambda_0}^{\lambda(s,r)} \sqrt{2s}\frac{M^2}{\mu}e^{\sqrt{2}r}d\mu ds+Mre^{\sqrt{2}r}\\
&\leq & M^2 e^{\sqrt{2}r}\int_{\varepsilon r}^{r-1}\sqrt{2s} \frac{\lambda(s,r)}{\lambda_0}ds+Mre^{\sqrt{2}r}.
\end{eqnarray*}
As $\lambda_0=\min\{1-\sqrt{2}/2,\,\sqrt{2/(r-1)}\}$, we deduce that $\mathbb{E}\bigg[e^{2\sigma_M-\sqrt{2}L(\sigma_M)};\Psi(r,\lambda)\bigg]\leq c_{13}M^2 r^2 e^{\sqrt{2}r}$, which completes the proof of (ii) in Lemma \ref{jointdis}. $\square$

Let us turn to prove the upper bound of $\Theta_k$.

\textit{Proof of (\ref{upperbound of main thm}).} For any $u\in\mathcal{N}(k)$ and $t>0$, we denote $A_u(t):=\{\tau_u>t\}$. Then $\{\Theta_k>t\}=\cup_{u\in\mathcal{N}(k)}A_u(t).$

For any $\theta\in \mathbb{Q}\cap(0,1)$, let
$$a_j:=\sqrt{2}-j\theta,\ \ b_j:=\sqrt{2}-(j-1)\theta, \text{ for } j=1,\dots, K:=K(\theta)=\lfloor\frac{\sqrt{2}}{\theta}\rfloor,$$
so that $0<a_j<\sqrt{2}$ for all $j\leq K$.

Let $I_k(a,b):=\{u\in\mathcal{N}(k); ak\leq X_u(k)\leq bk\}$ for $-\infty<a<b<\infty$.

Given the event $\Xi:=\Big\{-\sqrt{2}k\leq L(k)\leq R(k)\leq \sqrt{2}k\Big\}$, we can write 
$$\mathcal{N}(k)=I_k(-\theta,\sqrt{2})\cup\Big(\bigcup_{1\leq j\leq K}I_k(-b_j,-a_j)\Big),$$
\noindent so that
$$\{\Theta_k>t\}=\bigcup_{u\in\mathcal{N}(k)}A_u(t)=\Big( \bigcup_{u\in I_k(-\theta,\sqrt{2})}A_u(t)\Big)\cup\Big(\bigcup_{1\leq j\leq K}\bigcup_{u\in I_k(-b_j,-a_j)}A_u(t)\Big).$$
\noindent As a consequence,
\begin{equation}\label{grouping}
\mathbb{P}\bigg[\Big\{\Theta_k>t\Big\}\cap\Xi\bigg]\leq \mathbb{P}\Bigg[\bigcup_{u\in I_k(-\theta,\sqrt{2})}A_u(t)\cap \Xi\Bigg]+\sum_{j=1}^{K}\mathbb{P}\Bigg[\bigcup_{u\in I_k(-b_j,-a_j)}A_u(t)\cap\Xi\Bigg].
\end{equation}

We first estimate $\mathbb{P}\Big[\bigcup_{u\in I_k(-\theta,\sqrt{2})}A_u(t)\cap \Xi\Big]$.

For any particle $u\in\mathcal{N}(k)$, let $\sigma_M^u:=\inf\{s>0; N^u(s)=1+M\}$. Recall that $L^u(t)=\min\{X_v(t+k)-X_u(k); v\in\mathcal{N}(k+s), u<v\}$ for any $t>0$. By the branching property, conditioned on $\mathcal{F}_k$, $\{\sigma_M^u, L^u(\sigma^u_M)\}_{u\in\mathcal{N}(k)}$ are i.i.d. copies of $(\sigma_M, L(\sigma_M))$.

Similarly, We define $\Phi^u(r, \lambda):=\{\sigma^u_M>r-1\}\cup\{\varepsilon r\leq \sigma^u_M\leq r-1, L^u(\sigma^u_M)\leq-\lambda(r, \sigma^u_M)\sigma^u_M\}$ and $\Psi^u(r,\lambda):=\big\{\varepsilon r\leq \sigma_M^u \leq r-1, L^u(\sigma_M^u)\geq- \lambda(\sigma^u_M, r)\sigma_M^u \big\}$. One immediately observes that
\begin{multline}
\mathbb{P}\bigg[\bigcup_{u\in I_k(-\theta,\sqrt{2})}A_u(t)\cap \Xi\bigg]\leq \mathbb{P}\bigg[\bigcup_{u\in I_k(-\theta,\sqrt{2})}\Phi^u(r,\lambda)\bigg]\label{goodpoints}\\
+\mathbb{P}\bigg[\bigcup_{u\in I_k(-\theta,\sqrt{2})}A_u(t)\cap \bigg(\Phi^u(r,\lambda)\bigg)^c; \Xi\bigg].
\end{multline}
Conditioning on $\mathcal{F}_k$ yields that
\begin{equation}\label{lambdaa}
\mathbb{P}\Bigg[\bigcup_{u\in I_k(-\theta,\sqrt{2})}\Phi^u(r,\lambda)\Bigg]\leq\mathbb{E}\Bigg[\sum_{u\in I_k(-\theta,\sqrt{2})}1_{\Phi^u(r,\lambda)}\Bigg]=\mathbb{E}\Bigg[\sum_{u\in I_k(-\theta,\sqrt{2})}1\Bigg]\mathbb{P}\bigg[\Phi(r,\lambda)\bigg].
\end{equation}
Clearly,
\begin{equation}\label{largetype}
\mathbb{E}\Bigg[\sum_{u\in I_k(-\theta,\sqrt{2})}1\Bigg]\leq \mathbb{E}[N(k)]=e^k.
\end{equation}
It then follows from (\ref{phiupp}) that
\begin{equation}
\mathbb{P}\Big[\cup_{u\in I_k(-\theta,\sqrt{2})}\Phi^u(r,\lambda)\Big]\leq c_{12}M^2r e^{-r+k}.
\end{equation}
We now choose $r=k(1+\varepsilon)$ and set $\Lambda_0:=\mathbb{P}\Big[\cup_{u\in I_k(-\theta,\sqrt{2})}A_u(t)\cap \Big(\Phi^u(r,\lambda)\Big)^c; \Xi\Big]$. Then for all $k$ large enough, (\ref{goodpoints}) becomes
\begin{equation}\label{goodpointssimple}
\mathbb{P}\Bigg[\bigcup_{u\in I_k(-\theta,\sqrt{2})}A_u(t)\cap \Xi\Bigg]\leq c_{14}M^2 ke^{-\varepsilon k}+\Lambda_0.
\end{equation}

It remains to estimate $\Lambda_0$. Since $\Big(\Phi^u(r,\lambda)\Big)^c\subset \{\sigma^u_M<\varepsilon r\}\cup\Psi^u(r,\lambda)$, we write
\begin{equation}\label{lambdazero}
\Lambda_0=\mathbb{P}\Big[\bigcup_{u\in I_k(-\theta,\sqrt{2})}A_u(t)\cap \Big(\Phi^u(r,\lambda)\Big)^c; \Xi\Big]\leq \Lambda_1+\Lambda_2,
\end{equation}
where
\begin{eqnarray*}
\Lambda_1&:=&\mathbb{P}\bigg[\bigcup_{u\in I_k(-\theta,\sqrt{2})}A_u(t)\cap\{\sigma^u_M<\varepsilon r\}; \Xi\bigg]\\
\Lambda_2&:=&\mathbb{P}\bigg[\bigcup_{u\in I_k(-\theta,\sqrt{2})}A_u(t)\cap\Psi^u(r,\lambda); \Xi\bigg].
\end{eqnarray*}
For any particle $u\in I_k(-\theta,\sqrt{2})$ such that $\{\sigma_M^u<\varepsilon r\}$, for any $y\geq1$, we define
\begin{equation}
\widetilde{S}^u(y):=\min_{v\in\mathcal{N}^u(\varepsilon r)}T^v(y)+\varepsilon r.
\end{equation}
Recall that $A_u(t)=\{\tau_u>t\}$. By comparing $\widetilde{S}^u(y)$ with $t$, we obtain that for any $t_1\in (0, t)$,
\begin{eqnarray}\label{lambdaone}
\Lambda_1 &\leq&  \mathbb{P}\Big[\bigcup_{u\in I_k(-\theta,\sqrt{2})}\Big\{t_1\leq \widetilde{S}^u(y)\leq t<\tau_u\Big\};\, \Xi\Big]\\
&+&\mathbb{P}\Big[\bigcup_{u\in I_k(-\theta,\sqrt{2})}\Big\{\widetilde{S}^u(y)>t; \sigma^u_M<\varepsilon r\Big\}\Big]+\mathbb{P}\Big[\bigcup_{u\in I_k(-\theta,\sqrt{2})}\Big\{\widetilde{S}^u(y)<t_1\Big\}\Big]\nonumber\\
&=:&\Lambda_{1a}+\Lambda_{1b}+\Lambda_{1c}.\nonumber
\end{eqnarray}
For $\delta\in(0,1)$ and $\theta\in\mathbb{Q}\cap(0,1)$, we take $y=(\sqrt{2}+2)k(1+\theta)$, $t_1=e^k$, $t=e^{\sqrt{2}y(1+2\delta)}$. As $\{\sigma_M^u<\varepsilon r\}$ implies $\{N^u(\varepsilon r)> M\}$, $\mathbb{P}[\widetilde{S}^u(y)>t; \sigma^u_M<\varepsilon r\vert \mathcal{F}_{k+\varepsilon r}]$ is less than $\mathbb{P}[T(y)>t-\varepsilon r]^M$. Conditionally on $\mathcal{F}_{k+\varepsilon r}$, then by (\ref{largetype}), we get
\begin{equation}
\Lambda_{1b}\leq \mathbb{E}\Big[\sum_{u\in I_k(-\theta,\sqrt{2})}1\Big]\mathbb{P}\Big[T(y)>t-\varepsilon r\Big]^M\leq e^k\mathbb{P}\Big[T(y)>e^{\sqrt{2}y(1+\delta)}\Big]^M,
\end{equation}
since $r=(1+\varepsilon)k$ with $\epsilon\in(0,1/2)$.
By (\ref{largerpart}), $\Lambda_{1b}\leq e^k\times e^{-M\delta_1 y/3}$ for $k$ large enough. We take $M=\frac{6}{\delta_1}$ to ensure that $\Lambda_{1b}\leq e^{-k}$.

On the other hand, we observe that
\begin{equation}
\Lambda_{1c}\leq \mathbb{E}\bigg[\sum_{v\in\mathcal{N}(k+\varepsilon r)}1_{\big(T^v(y)< e^k\big)}\bigg]\leq e^{k+\varepsilon r}\mathbb{P}\bigg[T(y)<e^k\bigg].
\end{equation}
By (\ref{lesspart}), $\Lambda_{1c}\leq e^{k+\varepsilon r}c_5 y^2 e^{-\sqrt{2}y}e^k$. Thus $\Lambda_{1b}+\Lambda_{1c}\leq 2e^{-k}$ for sufficiently large $k$.

Set $\Xi_1:=\Big\{ R(k+\varepsilon r)\leq 2(k+\varepsilon r)\Big\}\cap\Xi$. Then,
\begin{equation}\label{lambdaonea}
\Lambda_{1a}\leq \mathbb{P}\Big[R(k+\varepsilon r)>2(k+\varepsilon r)\Big]+\mathbb{E}\Bigg[\sum_{u\in I_k (-\theta,\sqrt{2})}1_{(\tau_u>t\geq \widetilde{S}^u(y)\geq t_1)}; \Xi_1\Bigg].
\end{equation}
We define $\Lambda_{1rest}:=\mathbb{E}\Big[\sum_{u\in I_k (-\theta,\sqrt{2})}1_{(\tau_u>t\geq \tilde{S}^u(y)\geq t_1)}; \Xi_1\Big]$ for convenience. On the one hand, $\mathbb{P}[R(k+\varepsilon r)>2(k+\varepsilon r)]\leq e^{-k}$ because of (\ref{feynmankac}). On the other hand, we have
\begin{eqnarray*}
\Lambda_{1rest}&\leq &\mathbb{E}\Bigg[\sum_{u\in I_k(-\theta,\sqrt{2})} \int_{t_1}^{t} \mathbf{1}_{(\widetilde{S}^u(y)\in dr^\prime)}\mathbb{P}\bigg[\tau_u> r^\prime;\Xi_1\Big\vert \mathcal{F}_k,\mathcal{F}^u_\infty\bigg]\Bigg].
\end{eqnarray*}
Since $\{\widetilde{S}^u(y)=r^\prime\}\subset\{R^u(r^\prime)\geq L^u(\varepsilon r)+m(r^\prime-\varepsilon r)+y\}$, the event $\{\tau_u>r^\prime\}$ conditioned on $\{\widetilde{S}^u(y)=r^\prime\}$ implies $\cup_{w\in\mathcal{N}(k)\setminus\{ u\}}\Big\{X_w(k)+R^w(r^\prime)>X_u(k)+L^u(\varepsilon r)+m(r^\prime-\varepsilon r)+y\Big\}$. Further, this set is contained in $\cup_{w\in\mathcal{N}(k)\setminus\{u\}}\Big\{R^w(r^\prime)>m(r^\prime)+X_u(k)+L^u(\varepsilon r)+y-\sqrt{2}\varepsilon r-X_w(k)\Big\}$. As $\Xi_1$ guarantees that $\Big(X_u(k)+L^u(\varepsilon r)+y-\sqrt{2}\varepsilon r-X_w(k)\Big)_+\leq C_2 k$, the inequality (\ref{ineq}) can be applied to show that
\begin{eqnarray*}
&&\mathbb{E}\Bigg[\sum_{u\in I_k(-\theta,\sqrt{2})} \int_{t_1}^{t} \mathbf{1}_{(\widetilde{S}^u(y)\in dr^\prime)}\mathbb{P}\bigg[\tau_u> r^\prime;\Xi_1\Big\vert \mathcal{F}_k,\mathcal{F}^u_\infty\bigg]\Bigg]\\
& \leq & \mathbb{E}\Bigg[\sum_{u\in I_k(-\theta,\sqrt{2})} \int_{t_1}^{t} \mathbf{1}_{(\tilde{S}^u(y)\in dr^\prime)}\sum_{w\in\mathcal{N}(k)\setminus\{u\}} c_2(1+C_2 k)^2e^{-\sqrt{2}\Big(X_u(k)+L^u(\varepsilon r)+y-\sqrt{2}\varepsilon r-X_w(k)\Big)}\Bigg]\\
& \leq & c_{15}k^2e^{-\sqrt{2}y+2\varepsilon r}\mathbb{E}\Bigg[\sum_{u\in I_k(-\theta,\sqrt{2})}\sum_{w\in\mathcal{N}(k)\setminus\{u\}}e^{\sqrt{2}(X_w(k)-X_u(k))}\Bigg]\mathbb{E}\Bigg[e^{-\sqrt{2} L(\varepsilon r)}\Bigg],
\end{eqnarray*}
where the last inequality follows from the fact that for $u\in\mathcal{N}(k)$, $L^u(\varepsilon r)$ are independent of $\mathcal{F}_k$ and are independent copies of $L(\varepsilon r)$.

Whereas by the estimation of $\mathbb{E}\Big[\sum_{u\in\mathcal{N}_\delta(k)}\sum_{w\in\mathcal{N}(k)\setminus\{u\}}e^{\sqrt{2}X_w(k)}\Big]$ in Section 4,
\begin{eqnarray*}
&&\mathbb{E}\Bigg[\sum_{u\in I_k(-\theta,\sqrt{2})}\sum_{w\in\mathcal{N}(k)\setminus\{u\}}e^{\sqrt{2}(X_w(k)-X_u(k))}\Bigg]\\
&=&2\int_{0}^k e^{2k-s}ds\int_{\mathbb{R}}\mathbb{P}\Big[B_s\in dx\Big]\mathbb{E}\Big[e^{-\sqrt{2}B_k}; -\theta k\leq B_k\leq \sqrt{2}k\Big\vert B_s=x\Big]\mathbb{E}\Big[e^{\sqrt{2}B_k}\Big\vert B_s=x\Big]\\
&=&2\int_0^k e^{3k-2s}\mathbb{E}\Big[e^{\sqrt{2}B_s-\sqrt{2}B_k}; -\theta k\leq B_k\leq \sqrt{2}k\Big]ds.
\end{eqnarray*}
Because $\mathbb{E}\Big[e^{\sqrt{2}B_s-\sqrt{2}B_k}; -\theta k\leq B_k\leq \sqrt{2}k\Big]\leq e^{\sqrt{2}\theta k}\mathbb{E}\Big[e^{\sqrt{2}B_s}\Big]=e^{\sqrt{2}\theta k+s}$, we obtain that
\begin{equation}\label{manytotwo}
\mathbb{E}\Bigg[\sum_{u\in I_k(-\theta,\sqrt{2})}\sum_{w\in\mathcal{N}(k)\setminus\{u\}}e^{\sqrt{2}(X_w(k)-X_u(k))}\Bigg]\leq 2e^{3k+\sqrt{2}\theta k}.
\end{equation}
Besides, $\mathbb{E}\bigg[e^{-\sqrt{2} L(\varepsilon r)}\bigg]\leq \mathbb{E}\bigg[\sum_{v\in\mathcal{N}(\varepsilon r)}e^{-\sqrt{2}X_v(\varepsilon r)}\bigg]=e^{2\varepsilon r}$. As a result,
\begin{equation*}
\Lambda_{1rest}\leq c_{15}k^2e^{-\sqrt{2}y+2\varepsilon r}\times 2e^{3k+\sqrt{2}\theta k+2\varepsilon r}.
\end{equation*}
Recall that $0<\theta<1$, $y=(\sqrt{2}+2)k(1+\theta)$ and that $r=k(1+\varepsilon)$ with $k$ large enough so that $r>1/\varepsilon$. Hence, $\Lambda_{1rest}\leq c_{15} k^2 e^{(1-2\sqrt{2})k}$ for $\varepsilon\in(0,\frac{\theta}{3})$. Going back to (\ref{lambdaonea}), we get $\Lambda_{1a}\leq e^{-k}+c_{15} k^2 e^{(1-2\sqrt{2})k}$ for all $k$ sufficiently large.

Consequently, (\ref{lambdaone}) becomes
\begin{equation}\label{lambdaoneend}
\Lambda_1\leq e^{-k}+c_{15} k^2 e^{(1-2\sqrt{2})k}+2e^{-k}\leq c_{16}e^{-\varepsilon k}.
\end{equation}

It remains to estimate $\Lambda_2=\mathbb{P}\left[\bigcup_{u\in I_k(-\theta,\sqrt{2})}A_u(t)\cap\Psi^u(r,\lambda);\Xi\right]$ where $\Psi^u(r,\lambda)=\big\{\varepsilon r\leq \sigma_M^u \leq r-1, L^u(\sigma_M^u)\geq- \lambda(\sigma^u_M, r)\sigma_M^u \big\}$. For any particle  $u\in\mathcal{N}(k)$ satisfying $\Psi^u(r, \lambda)$, define
\begin{equation}\label{hattime}
\widehat{S}^u(y):=\min_{v\in\mathcal{N}^u(\sigma^u_M)}T^{v}(y)+\sigma^u_M,\text{  for any }y>0.
\end{equation}
Comparing $\widehat{S}^u(y)$ with $t$ yields that
\begin{multline}\label{lambdatwo}
\Lambda_2 \leq  \mathbb{P}\bigg[\bigcup_{u\in I_k(-\theta,\sqrt{2})}\Big\{t_1\leq \widehat{S}^u(y)\leq t<\tau_u\Big\}\cap\Psi^u(r,\lambda); \Xi\bigg]\\
+\mathbb{P}\bigg[\bigcup_{u\in I_k(-\theta,\sqrt{2})}\Big\{\widehat{S}^u(y)>t\Big\}\cap\Psi^u(r,\lambda)\bigg]+\mathbb{P}\bigg[\bigcup_{u\in I_k(-\theta,\sqrt{2})}\Big\{\widehat{S}^u(y)<t_1\Big\}\cap\Psi^u(r,\lambda)\bigg]\\
=:\Lambda_{2a}+\Lambda_{2b}+\Lambda_{2c}.
\end{multline}
According to the definition of $\widehat{S}^u(y)$, one sees that $\mathbb{P}[\widehat{S}^u(y)>t;\; \Psi^u(r,\lambda)]\leq\mathbb{P}[T(y)>t-r]^M$ and that $\mathbb{P}[\widehat{S}^u(y)<t_1;\; \Psi^u(r,\lambda)]\leq 2M\mathbb{P}[T(y)<t_1]$. Recall that $r=(1+\varepsilon)k$, $M=6/\delta_1$, $y=(\sqrt{2}+2)k(1+\theta)$, $t_1=e^k$ and $t=e^{\sqrt{2}y(1+2\delta)}$. For any $-\infty<a<b<\infty$, by (\ref{largerpart}),
\begin{equation}\label{lambda2b}
 \mathbb{E}\Bigg[\sum_{u\in I_k(a,b)}1_{\{\widehat{S}^u(y)>t\}\cap\Psi^u(r,\lambda)}\Bigg]\leq e^k\times\mathbb{P}[T(y)>e^{\sqrt{2}y(1+\delta)}]^M\leq e^{-k}.
\end{equation}
Meanwhile, by (\ref{lesspart}),
\begin{equation}\label{lambda2c}
\mathbb{E}\Bigg[\sum_{u\in I_k(a,b)}1_{\{\widehat{S}^u(y)<t_1\}\cap\Psi^u(r,\lambda)}\Bigg]\leq e^k\times2M\mathbb{P}[T(y)<t_1]\leq 2c_5Me^{-2k}.
\end{equation}
Hence, taking $a=-\theta$ and $b=\sqrt{2}$ implies that $\Lambda_{2b}+\Lambda_{2c}\leq e^{-k}+2c_5Me^{-2k}$.
Let $\Xi_2:=\{\max_{0\leq r_0\leq r} R(k+r_0)\leq 6k\}\cap \Xi$. We get
\begin{equation}\label{lambda2a}
\Lambda_{2a}\leq \mathbb{P}\Big[\max_{0\leq r_0\leq r} R(k+r_0)>6k\Big]+\mathbb{E}\Bigg[\sum_{u\in I_k(-\theta,\sqrt{2})}1_{\{t_1\leq \widehat{S}^u(y)\leq t<\tau_u\}\cap\Psi^u(r,\lambda)};\, \Xi_2\Bigg].
\end{equation}
By the many-to-one lemma, for $k$ large enough,
\begin{equation}\label{xitwo}
\mathbb{P}\Big[\max_{0\leq r_0\leq r} R(k+r_0)>6k\Big]\leq e^{-k}.
\end{equation}
For the second term on the right-hand side of (\ref{lambda2a}), we need to recount the arguments to estimate $\Lambda_{1rest}$. Let
\begin{equation}
\Lambda_{2rest}:=\mathbb{E}\Bigg[\sum_{u\in I_k(-\theta,\sqrt{2})}1_{\{t_1\leq \widehat{S}^u(y)\leq t<\tau_u\}\cap\Psi^u(r,\lambda)};\, \Xi_2\Bigg].
\end{equation}
It immediately follows that
\begin{equation}
\Lambda_{2rest}\leq \mathbb{E}\left[\sum_{u\in I_k(-\theta,\sqrt{2})}\int_{t_1}^t 1_{(\widehat{S}^u(y)\in dr^\prime)}1_{\Psi^u(r,\lambda)}\times\mathbb{P}\bigg[\{\tau_u>r^\prime\}\cap\Xi_2\Big\vert \mathcal{F}_k, \mathcal{F}^u_\infty\bigg] \right].
\end{equation}
Comparing $\tau_u$ with $\widehat{S}^u(y)$ tells that
\begin{multline*}
\Lambda_{2rest}\leq  \mathbb{E}\Bigg[\sum_{u\in I_k(-\theta,\sqrt{2})}\int_{t_1}^t 1_{\big(\widehat{S}^u(y)\in dr^\prime\big)}1_{\Psi^u(r,\lambda)}\times\\
\mathbb{E}\bigg[1_{\Xi_2}\sum_{w\in\mathcal{N}(k)\setminus\{u\}}1_{\big(R^w(r^\prime)>m(r^\prime)-\sqrt{2}\sigma^u_M+y+X_u(k)-X_w(k)+L^u(\sigma^u_M)\big)}\Big\vert\mathcal{F}_k, \mathcal{F}^u_\infty\bigg] \Bigg].
\end{multline*}
On the event $\Xi_2\cap\Psi^u(r,\lambda)$, we have $1+\Big(-\sqrt{2}\sigma^u_M+y+X_u(k)-X_w(k)+L^u(\sigma^u_M)\Big)_+\leq C_3 k$. Applying the inequality (\ref{ineq}) for $R^w(r^\prime)$ yields that
\begin{eqnarray*}
\Lambda_{2rest}& \leq & \mathbb{E}\Bigg[\sum_{u\in I_k(-\theta,\sqrt{2})} 1_{\Psi^u(r,\lambda)}\times\sum_{w\in\mathcal{N}(k)\setminus\{u\}}C_3^2 k^2e^{-\sqrt{2}y}e^{\sqrt{2}(X_w(k)-X_u(k))}e^{2\sigma_M^u-\sqrt{2}L^u(\sigma_M^u)}\Bigg]\\
&=& C_3^2 k^2 e^{-\sqrt{2}y}\mathbb{E}\Bigg[\sum_{u\in I_k(-\theta,\sqrt{2})}\sum_{w\in\mathcal{N}(k)\setminus\{u\}}e^{\sqrt{2}(X_w(k)-X_u(k))}\Bigg]\mathbb{E}\Bigg[e^{2\sigma_M-\sqrt{2}L(\sigma_M)};\Psi(r,\lambda)\Bigg],
\end{eqnarray*}
by the fact that $\{\sigma_M^u, L^u(\sigma_M^u)\}$ are i.i.d. and independent of  $\mathcal{F}_k$.
Recall that $r=(1+\varepsilon)k$ with $\varepsilon\in(0,\theta/3)$. It then follows from (\ref{manytotwo}) and (\ref{psiupp}) that
\begin{equation}
\Lambda_{2rest}\leq C_3^2 k^2 e^{-(2+2\sqrt{2})k(1+\theta)}\times2e^{3k+\sqrt{2}\theta k}\times c_{13}M^2 r^2e^{\sqrt{2}r}\leq c_{17} k^4 M^2 e^{-(\sqrt{2}-1) k}.
\end{equation}
Consequently, $\Lambda_{2a}\leq e^{-k}+c_{17} k^4 M^2 e^{-(\sqrt{2}-1) k}$.

Therefore, for $k$ sufficiently large, the inequality (\ref{lambdatwo}) becomes
\begin{equation}\label{lambdatwoend}
\Lambda_2\leq e^{-k}+2c_5Me^{-2k}+e^{-k}+c_{17} k^4 M^2 e^{-(\sqrt{2}-1) k}\leq c_{18}M^2e^{-\varepsilon k}.
\end{equation}
Combined with (\ref{lambdaoneend}), $\Lambda_0\leq c_{16}e^{-\varepsilon k}+c_{18}M^2e^{-\varepsilon k}$. Going back to (\ref{goodpointssimple}), we conclude that
\begin{equation}\label{goodpointsconclusion}
\mathbb{P}\Bigg[\bigcup_{u\in I_k(-\theta,\sqrt{2})}A_u(t)\cap \Xi\Bigg]\leq c_{19}M^2e^{-\varepsilon k/2}.
\end{equation}

To complete the proof, we still need to evaluate $\mathbb{P}\Big[\cup_{u\in I_k(-b_j,-a_j)}A_u(t)\cap\Xi\Big]$. Recall that for any particle $u\in\mathcal{N}(k)$, $\sigma_M^u=\inf\{s>0; N^u(s)=1+M\}$ and
$$\Phi^u(r, \lambda)=\{\sigma^u_M>r-1\}\cup\{\varepsilon r\leq \sigma^u_M\leq r-1, L^u(\sigma^u_M)\leq-\lambda(\sigma^u_M,r)\sigma^u_M\},$$
for any $r>1/\varepsilon$ and $\lambda(s, r)=\sqrt{2(\frac{r}{s}-1)}$ with $0<s<r$. Clearly,
\begin{equation}
\mathbb{P}\Bigg[\bigcup_{u\in I_k(-b_j,-a_j)}A_u(t);\Xi\Bigg]\leq \mathbb{P}\Bigg[\bigcup_{u\in I_k(-b_j,-a_j)}\Phi^u(r, \lambda)\Bigg]+\mathbb{P}\Bigg[\bigcup_{u\in I_k(-b_j,-a_j)}A_u(t)\cap\Big(\Phi^u(r,\lambda)\Big)^c;\Xi\Bigg].\nonumber
\end{equation}
On the one hand,
\begin{equation}\label{phisum}
\mathbb{P}\Bigg[\bigcup_{u\in I_k(-b_j,-a_j)}\Phi^u(r,\lambda)\Bigg]\leq\mathbb{E}\Bigg[\sum_{u\in I_k(-b_j,-a_j)}1\Bigg]\mathbb{P}\bigg[\Phi(r,\lambda)\bigg].
\end{equation}
We now take $r=k(1-\frac{a_j^2}{2})(1+\varepsilon)$ with $\varepsilon>0$ small so that $r\leq 2k$. Recall that $a_j=\sqrt{2}-j\theta$ for $1\leq j\leq \lfloor \frac{\sqrt{2}}{\theta}\rfloor$, with $\theta\in\mathbb{Q}\cap (0,1)$. Then note that each $a_j$ is strictly positive. Thus, by the many-to-one lemma and (\ref{phiupp}), (\ref{phisum}) becomes that
\begin{eqnarray*}
\mathbb{P}\Bigg[\bigcup_{u\in I_k(-b_j,-a_j)}\Phi^u(r,\lambda)\Bigg]&\leq& e^k \mathbb{P}[-b_j k\leq B_k\leq -a_j k]\times c_{12}M^2re^{-r}\\
&\leq&  c_{12}M^2re^{-r+k}\mathbb{P}[B_k\leq -a_j k]\\
&\leq& c_{12}M^2 (2k) e^{-k(1-\frac{a_j^2}{2})(1+\varepsilon)+k}\Big(\frac{\sqrt{k}}{a_j k}e^{-\frac{a_j^2}{2}k}\Big),
\end{eqnarray*}
where the last inequality follows from (\ref{evaluationofBM}). $\mathbb{P}\Bigg[\bigcup_{u\in I_k(-b_j,-a_j)}\Phi^u(r,\lambda)\Bigg]$ is hence bounded by $c_{12}\frac{2M^2}{a_j}\sqrt{k}e^{-\sqrt{2}\theta\varepsilon k/2}$.

On the other hand, recalling that $\Psi^u(r,\lambda)=\{\varepsilon r\leq \sigma^u_M\leq r-1, L^u(\sigma^u_M)\geq-\lambda(\sigma^u_M,r)\sigma^u_M\}$, we deduce that
\begin{equation}\label{lambdazeroprime}
\mathbb{P}\Bigg[\bigcup_{u\in I_k(-b_j,-a_j)}A_u(t)\cap\Big(\Phi^u(r,\lambda)\Big)^c;\Xi\Bigg]\leq\Lambda_1^\prime+\Lambda_2^\prime,
\end{equation}
where
\begin{eqnarray*}
\Lambda_1^\prime &:=&\mathbb{P}\Bigg[\bigcup_{u\in I_k(-b_j, -a_j)}A_u(t)\cap\bigg\{\sigma_M^u<\varepsilon r\bigg\};\Xi\Bigg],\\
\Lambda_2^\prime &:=&\mathbb{P}\Bigg[\bigcup_{u\in I_k(-b_j, -a_j)}A_u(t)\cap\Psi^u(r,\lambda);\Xi\Bigg].
\end{eqnarray*}
Furthermore, by an argument similar to the one used in estimating $\Lambda_1$, we have $\Lambda_1^\prime\leq c_{20}e^{-\varepsilon k}$. Thus,
\begin{equation}\label{badpoints}
\mathbb{P}\Bigg[\bigcup_{u\in I_k(-b_j,-a_j)}A_u(t);\Xi\Bigg]\leq c_{12}\frac{2M^2}{a_j}\sqrt{k}e^{-\sqrt{2}\theta\varepsilon k/2}+c_{20}e^{-\varepsilon k}+\Lambda_2^\prime.
\end{equation}
It remains to bound $\Lambda_2^\prime$. Recall that $\widehat{S}^u(y)=\min_{v\in\mathcal{N}^u(\sigma^u_M)}T^{v}(y)+\sigma^u_M$ with $y=(\sqrt{2}+2)k(1+\theta)$. We observe that
\begin{eqnarray}
\Lambda_2^\prime& \leq & \mathbb{P}\Bigg[\bigcup_{u\in I_k(-b_j,-a_j)}\{\widehat{S}^u(y)>t\}\cap\Psi^u(r,\lambda)\Bigg]+\mathbb{P}\Bigg[\bigcup_{u\in I_k(-b_j,-a_j)}\{\widehat{S}^u(y)<t_1\}\cap\Psi^u(r,\lambda)\Bigg]\nonumber\\
& + & \mathbb{P}\Big[\max_{0\leq r_0\leq r} R(k+r_0)>6k\Big]+\mathbb{P}\Bigg[\bigcup_{u\in I_k(-b_j,-a_j)}\Big\{t_1\leq \widehat{S}^u(y)\leq t<\tau_u\Big\}\cap\Psi^u(r,\lambda);\Xi_2\Bigg].\nonumber
\end{eqnarray}
In view of (\ref{lambda2b}), (\ref{lambda2c}) and (\ref{xitwo}),
\begin{equation}
\Lambda_2^\prime\leq e^{-k}+2c_{5}M e^{-2k}+e^{-k}+\mathbb{E}\bigg[\sum_{u\in I_k(-b_j,-a_j)}1_{\{t_1\leq \widehat{S}^u(y)\leq t<\tau_u\}\cap\Psi^u(r,\lambda)};\Xi_2\bigg].
\end{equation}
We define $\Lambda_{2rest}^\prime:=\mathbb{E}\Big[\sum_{u\in I_k(-b_j,-a_j)}1_{\{t_1\leq \widehat{S}^u(y)\leq t<\tau_u\}\cap\Psi^u(r,\lambda)};\;\Xi_2\Big]$. Thus applying the analogous arguments to the estimation of $\Lambda_{2rest}$ gives that
\begin{multline}\label{lambda-2rest}
\Lambda_{2rest}^\prime\leq c_{21} k^2 e^{-\sqrt{2}y}\mathbb{E}\Bigg[\sum_{u\in I_k(-b_j,-a_j)}\sum_{w\in\mathcal{N}(k)\setminus\{u\}}e^{\sqrt{2}(X_w(k)-X_u(k))}\Bigg]\\
\times\mathbb{E}\Bigg[e^{2\sigma_M-\sqrt{2}L(\sigma_M)};\Psi^u(r,\lambda)\Bigg].
\end{multline}
Once again, by means of integrating with respect to the last time at which the most recent common ancestor of $u$ and $v$ was alive, $\mathbb{E}\Big[\sum_{u\in I_k(-b_j,-a_j)}\sum_{w\in\mathcal{N}(k)\setminus\{u\}}e^{\sqrt{2}(X_w(k)-X_u(k))}\Big]$ equals
\begin{eqnarray*}
&&2\int_0^k e^{2k-s}ds\int_{\mathbb{R}}\mathbb{P}\Big[B_s\in dx\Big]\mathbb{E}\Big[e^{-\sqrt{2}B_k}; -b_j k\leq B_k\leq -a_j k\Big\vert B_s=x\Big]\mathbb{E}\Big[e^{\sqrt{2}B_k}\Big\vert B_s=x\Big]\\
&=&2\int_0^k e^{3k-2s}\mathbb{E}\Big[e^{\sqrt{2}B_s-\sqrt{2}B_k}; -b_j k\leq B_k\leq -a_j k\Big]ds\\
&=&2\int_0^k e^{3k-2s}ds\int_{-b_jk}^{-a_jk}e^{-\sqrt{2}x}\mathbb{P}\Big[B_k\in dx\Big]\mathbb{E}\Big[e^{\sqrt{2}B_s}\Big\vert B_k=x\Big].
\end{eqnarray*}
Let $(b_s(x);\, 0\leq s\leq k)$ denote a Brownian bridge from $0$ to $x$ of length $k$. Then $\mathbb{E}\Big[e^{\sqrt{2}B_s}\Big\vert B_k=x\Big]$ equals $\mathbb{E}\Big[e^{\sqrt{2}b_s(x)}\Big]$, which turns out to be $\exp(s(k-s+\sqrt{2}x)/k)$. Note that $a_j>0$ and that $b_j=a_j+\theta$. This gives that
\begin{eqnarray*}
&&\mathbb{E}\Big[\sum_{u\in I_k(-b_j,-a_j)}\sum_{w\in\mathcal{N}(k)\setminus\{u\}}e^{\sqrt{2}(X_w(k)-X_u(k))}\Big]\\
&=& 2\int_0^k e^{3k-2s+s(k-s)/k}\int_{-b_jk}^{-a_jk}\frac{1}{\sqrt{2\pi k}}\exp\Big(-\frac{x^2}{2k}-\sqrt{2}x+\sqrt{2}\frac{s}{k}x\Big)dx ds\\
&\leq & 2\int_0^k e^{3k-2s+s(k-s)/k} \frac{\theta k}{\sqrt{2\pi k}}\exp\Big(-\frac{a_j^2}{2}k+\sqrt{2}b_jk-\sqrt{2}a_js\Big) ds,
\end{eqnarray*}
which is bounded by $\sqrt{k}\exp\Big(3k-\frac{a_j^2}{2}k+\sqrt{2}b_jk\Big)$ since $\int_0^k \exp\Big(-2s+s(k-s)/k-\sqrt{2}a_js\Big)ds$ is less than $1$. One hence sees that $\mathbb{E}\Big[\sum_{u\in I_k(-b_j,-a_j)}\sum_{w\in\mathcal{N}(k)\setminus\{u\}}e^{\sqrt{2}(X_w(k)-X_u(k))}\Big]$ is bounded by $\sqrt{k}e^{2k+(1-\frac{a_j^2}{2}+\sqrt{2}a_j)k+\sqrt{2}\theta k}$. Going back to (\ref{lambda-2rest}) and applying (\ref{psiupp}),
\begin{eqnarray*}
\Lambda_{2rest}^\prime& \leq & c_{21} k^2 e^{-(2+2\sqrt{2})k(1+\theta)}\sqrt{k}e^{2k+(1-\frac{a_j^2}{2}+\sqrt{2}a_j)k+\sqrt{2}\theta k} c_{13}M^2 r^2e^{\sqrt{2}r}\\
& \leq & c_{22}k^{9/2}M^2e^{-2\theta k}\exp\bigg\{k\Big[(\sqrt{2}+1)(1-a_j^2/2)+\sqrt{2}a_j-2\sqrt{2}\Big]\bigg\},
\end{eqnarray*}
as $r=k(1-\frac{a_j^2}{2})(1+\varepsilon)$ with $\varepsilon\in(0,\frac{\theta}{3})$. Observe that
\begin{equation}\label{optimal}
(\sqrt{2}+1)(1-a_j^2/2)+\sqrt{2}a_j-2\sqrt{2}=-\frac{\sqrt{2}+1}{2}\Big(a_j-(2-\sqrt{2})\Big)^2\leq 0.
\end{equation}
We get $\Lambda_{2rest}^\prime\leq c_{22}k^{9/2}M^2e^{-2\theta k}$, and thus for all $k$ sufficiently large,
\begin{equation}
\Lambda_2^\prime\leq c_{23}M^2 e^{-\theta k}.
\end{equation}
Consequently, by (\ref{badpoints}),
\begin{equation*}
\mathbb{P}\bigg[\cup_{u\in I_k(-b_j,-a_j)}A_u(t);\Xi\bigg]\leq c_{12}\frac{2M^2}{a_j}\sqrt{k}e^{-\sqrt{2}\theta \varepsilon k/2}+c_{20}e^{-\varepsilon k}+c_{23}M^2e^{-\theta k}.
\end{equation*}
Summing over $j\in\{1,\cdots,K=\lfloor \frac{\sqrt{2}}{\theta}\rfloor\}$ implies that
\begin{equation}\label{badpointsconclusion}
\sum_{j=1}^K \mathbb{P}\bigg[\cup_{u\in I_k(-b_j,-a_j)}A_u(t);\Xi\bigg]\leq C(\theta) M^2 e^{-\varepsilon\theta k/2},
\end{equation}
where $C(\theta)$ is a positive constant associated with $\theta$ (but independent of $k$, $\delta$ and $M$) and $k$ is large enough.

Going back to (\ref{grouping}), we combine (\ref{goodpointsconclusion}) and (\ref{badpointsconclusion}) to say that
\begin{eqnarray*}
\mathbb{P}\Bigg[\bigg\{\tau(k)>t\bigg\}\cap\bigg\{-\sqrt{2}k\leq L(k)\leq R(k)\leq \sqrt{2}k\bigg\}\Bigg]&\leq &c_{19} M^2 e^{-\varepsilon k/2}+C(\theta) M^2 e^{-\varepsilon\theta k/2}\\
&\leq &\frac{ C_1(\theta)}{\delta^2}e^{-\varepsilon\theta k/2},
\end{eqnarray*}
where $\theta\in\mathbb{Q}\cap(0,1)$, $\varepsilon\in(0,\frac{\theta}{3})$, $\delta>0$ and $t=\exp[k(2+2\sqrt{2})(1+\theta)(1+2\delta)]$ and $k$ is sufficiently large.

According to the Borel-Cantelli Lemma, we conclude that for any $\theta\in\mathbb{Q}\cap(0,1)$ and any $\delta>0$,
$$\limsup_{k\rightarrow\infty}\frac{\log\Theta_k}{k}\leq (2+2\sqrt{2})(1+\theta)(1+2\delta),\ \quad \textrm{almost surely.}$$
This implies the upper bound in Theorem \ref{mainconclusion}.  $\square$

\textbf{Acknowledgements}

I would like to thank Zhan Shi for advice and help.

\bibliographystyle{plainnat}

\bigskip
\bigskip
\end{document}